\documentclass[12pt,oneside]{article}
\newcommand{\Path }{}
\newcommand{\figs}{}
\usepackage{\Path Paper_style}
\usepackage{\Path Keywords}
\usepackage{\Path Environments}
\usepackage{tikz-cd}
\usepackage{tabularx} \newcolumntype{L}{>{\raggedright\arraybackslash}X}
\newcommand{\rkhs}{\mathcal{H}}

\newcommand{\Smooth}{\mathcal{G}}
\newcommand{\noise}{\eta}

\begin{document}
	\title{\vspace{-1cm} A probabilistic approach to drift estimation from stochastic data }
	\author{Suddhasattwa Das\footnotemark[1]}
	\footnotetext[1]{Department of Mathematics and Statistics, Texas Tech University, Texas, USA}
	\date{\today}
	\maketitle
	\begin{abstract} Timeseries generated from a dynamical source can often be modeled as sample paths of a stochastic differential equation (SDE). The timeseries thus reflects the motion of a particle which flows along the direction provided by a drift / vector field, and is simultaneously scattered by the effect of white noise. The resulting motion can only be described as a random process instead of a solution curve. Due to the non-deterministic nature of this motion, the task of determining the drift from data is quite challenging, since the data does not directly represent the directional information of the flow. This paper describes an interpretation of a drift as a conditional expectation, which makes its estimation feasible via kernel-integral methods. In addition, some techniques are proposed to overcome the challenge of dimensionality if the SDE's carry some structure enabling sparsity. The technique is shown to be convergent, consistent and permits a wide choice of kernels.
	\end{abstract}
	
	
	\section{Introduction} \label{sec:intro}
	
	The dynamics present in several physical phenomenon are governed by both deterministic laws of motion, as well as stochastic driving terms. Stochastic differential equations (SDEs) provide a common mathematical description for these phenomenon. A general $d$-dimensional SDE takes the form
	\begin{equation} \label{eqn:sde}
		dX(t) = V \paran{ X(t) } + G \paran{ X(t) } dW^t ,
	\end{equation}
	in which the term $V$ known as the \emph{drift}, is a function of the current state. It plays the role of a vector field over the state space $\real^d$, providing a deterministic direction to the flow at all locations. The term $G$ is known as the \emph{diffusion} term. It relays the effect of an external white noise onto the motion. 
	
	\begin{figure}[!t]\center	         
		\includegraphics[width=0.9\linewidth, height=0.3\textheight, keepaspectratio]{\figs 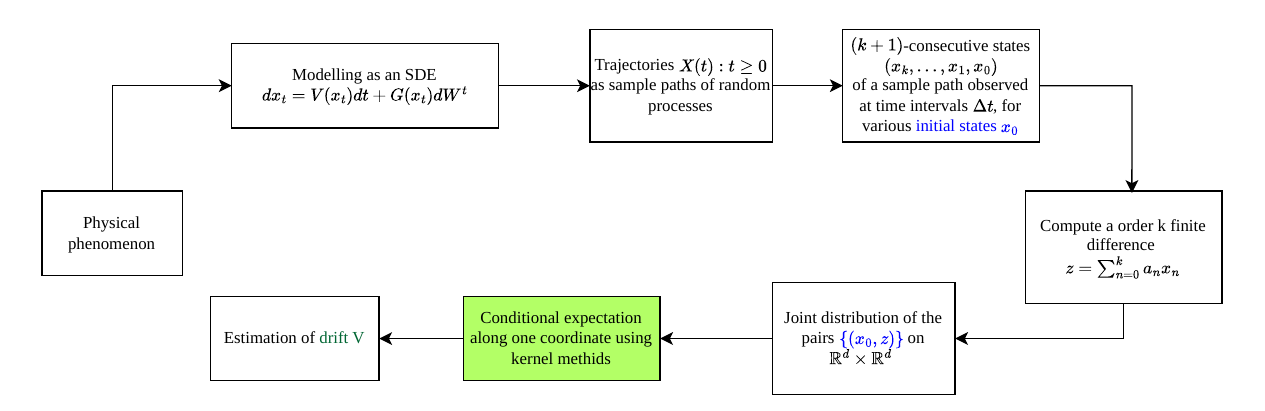}
		\caption{Outline of the theory. Many dynamical processes in nature can be modeled as an SDE, which has a drift and diffusion to it. Unlike deterministic systems, the consecutive points on a trajectory or not related deterministically, but have a distribution. The flowchart shows how we can collect $(k+1)$-length snapshots of this process, and apply purely data-driven algorithms to estimate the drift. No presumption is needed on a prior distribution on the coordinates, or on a parametric form for the SDE terms. The article presents the details of these methods, their theoretical footings in Probability theory, and a rigorous demonstration that they provide a convergent, unbiased estimate.}
		\label{fig:outline1}
	\end{figure}
	
	This article focuses on the data-driven discovery of SDEs. A convergent and consistent technique is presented for estimating the drift term $V$ in SDEs, from data collected from observations of the system at regular sampling intervals of $\Delta t$. The drift is a vector field and thus a deterministic function of the state. However due to the stochasticity originating from white noise the state at time $t$ of an SDE-driven system is not a deterministic function of its current state. Rather it has a probability distribution conditioned on the current state. As a result the solution $X(t,x_0)$ when expressed a function of time and initial state $x_0$ is not a deterministic curve but a Markov process in which almost every sample path is continuous but non-differentiable \cite[e.g.]{Arnold1974stoch, philipp2023error}. This stochasticity and lack of differentiability makes the task of drift estimation challenging, as the drift is essentially a differential operator.  
	
	The Weiner process that drives an SDE of the form \eqref{eqn:sde} can also be interpreted as a probability space, whose domain is the \emph{path-space} $\Omega_{\text{path}}$ of all possible random Brownian walks. The space $\Omega_{\text{path}}$ provides a concrete way to identity the source of stochasticity in the solutions of the SDE. Every random point from $\Omega_{\text{path}}$ corresponds to a particular sample of a Brownian path, and thus corresponds to a unique realization of \eqref{eqn:sde}. This path space formalism of stochastic processes, championed by Doob \cite{Doob1953book} and Arnold \cite{Arnold1974stoch} clarifies the source of stochasticity in the study of SDEs. The path space $\Omega_{\text{path}}$ itself is usually left undetermined, and the SDE is studied via the distributions it induces on $\real^d$. 
	
\begin{figure}[!t]\center
		\includegraphics[width=0.3\linewidth, height=0.3\textheight, keepaspectratio]{\figs 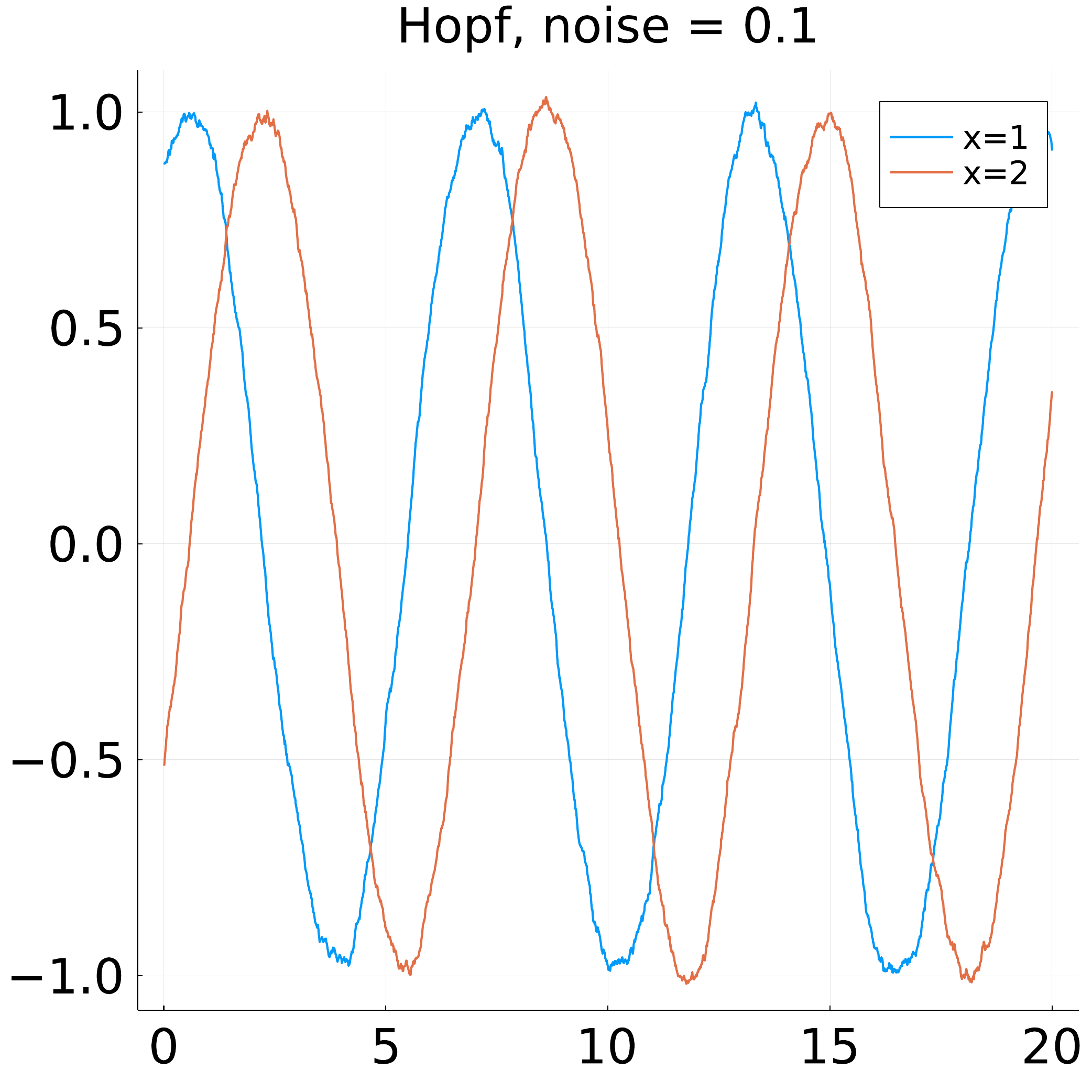}
		\includegraphics[width=0.3\linewidth, height=0.3\textheight, keepaspectratio]{\figs 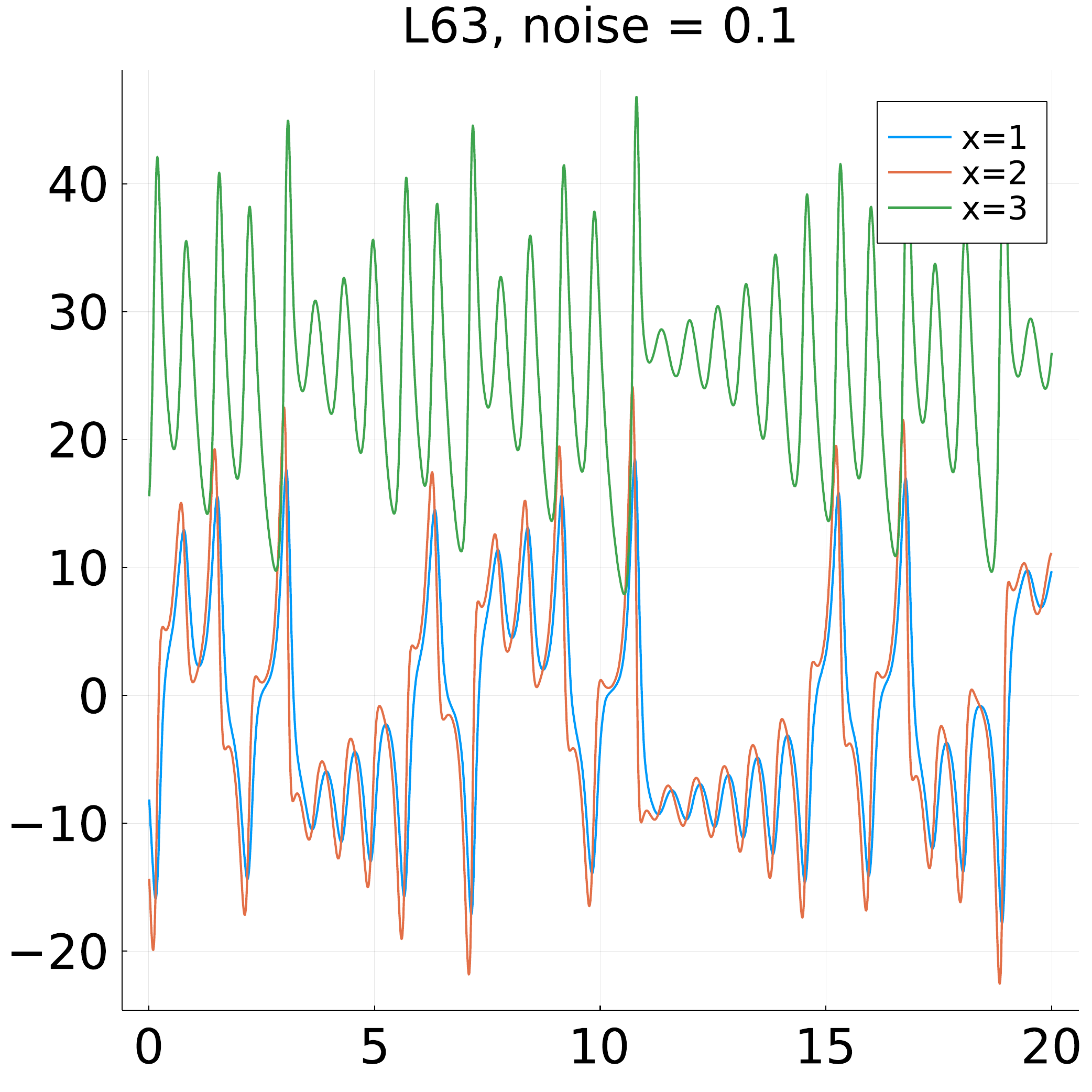}
		\includegraphics[width=0.3\linewidth, height=0.3\textheight, keepaspectratio]{\figs 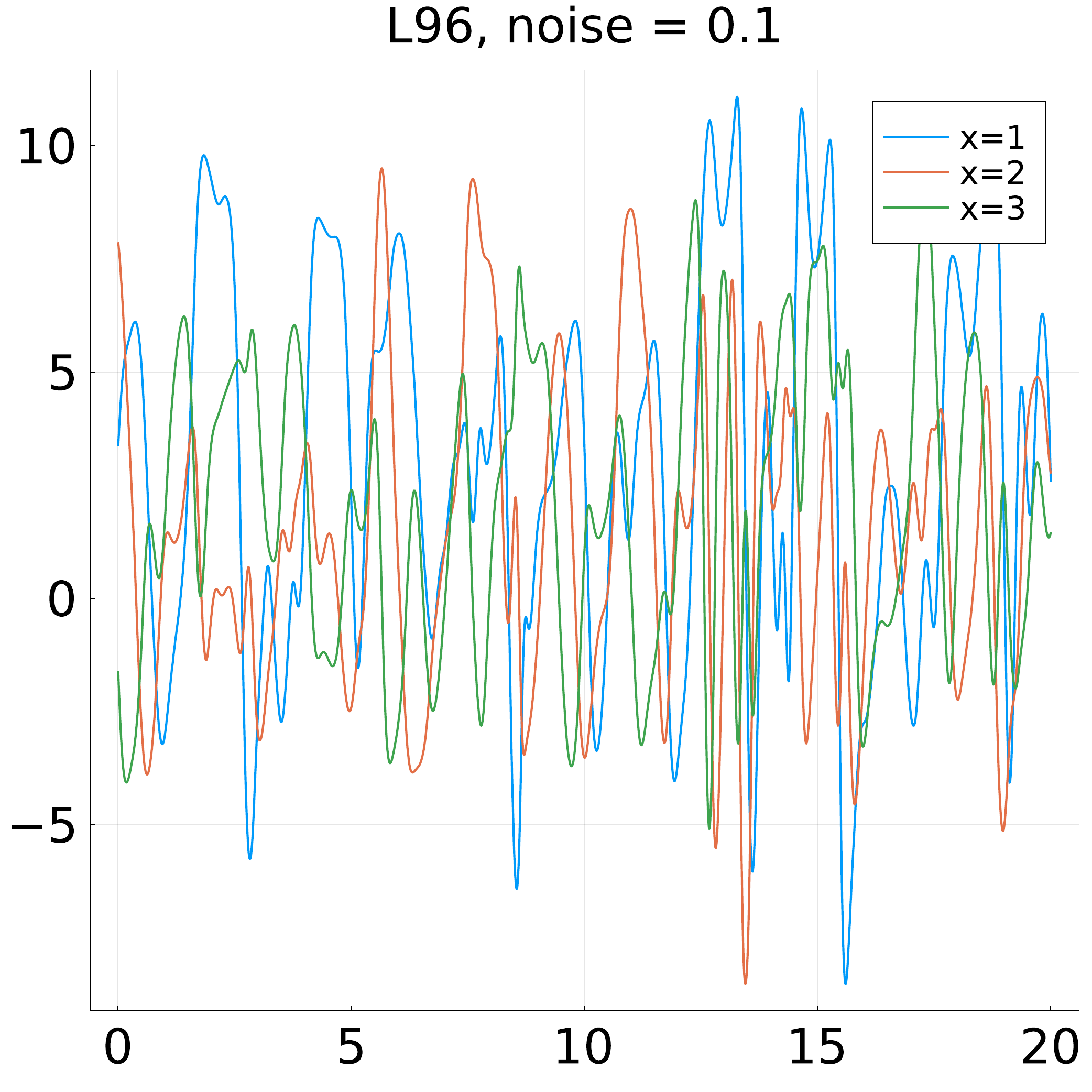}
		\caption{Orbits from SDEs. The three plots show sample paths of SDEs, in which the drifts are respectively the Hopf oscillator \eqref{eqn:Hopf}, Lorenz-63 system \eqref{eqn:vec:l63} and the Lorenz 96 system \eqref{eqn:L96}. The SDEs have a diagonal diffusion term given by \eqref{eqn:dffsn:expt}. The parameters of the experiments are provided in Table \ref{tab:param1}. The legends indicate the index of the coordinates being plotted. The paper presents a numerical method for extracting the drift from such timeseries which are generated by SDEs. The results of the reconstruction have been displayed in Figures \ref{fig:l63}, \ref{fig:Hopf}, and \ref{fig:L96_N5} respectively. }
		\label{fig:orbits}
	\end{figure}

	Drift estimation of SDE timeseries require a treatment completely different from those used for numerical differentiation. Figure \ref{fig:outline1} provides an outline of the mathematical problem, and the proposed method of converting it into a task of conditional expectation. In the usual probabilistic approach to SDEs, the initial state $x_0$ to the SDE \eqref{eqn:sde} is not fixed but is itself a random variable $X_0$. This stochasticity in the location of $x_0$ is propagated into the future states in combination with the effect of white-noise. As a result the state at any time $t$ is a random variable $X_t$. This 1-parameter family $\SetDef{X_t}{t\geq 0}$ of random variables is a Markov process. Its key feature is the useful \emph{memoryless property} -- for any three time instants $t_1 < t_2 < t_3$, the distribution of $X_{t_3}$ can be determined from the value of $X_{t_2}$ alone, and may be treated independent of $X_{t_1}$. Figure \ref{fig:orbits} presents some examples of SDEs simulated using the Euler Mauryama method. We make the following assumption  :
	
	\begin{Assumption} \label{A:1}
		The SDE \eqref{eqn:sde} has a stationary probability measure $\mu_{stat}$.
	\end{Assumption}
	
	Stationarity of $\mu_{stat}$ means that if $X_0$ is distributed according to $\mu_{stat}$, then so will be $X_t$ for every $t\geq 0$. Stationary measures are the stochastic analog of \emph{ergodic measures} in deterministic dynamical systems, and is one of the pillars of data-driven studies of stochastic systems. Stationarity is empirically observed in most systems, although it has been theoretically proven only under certain sets of conditions \cite[e.g.]{HuangEtAl2015integral, harris1956existence, HuangEtAl2018concn}.
	Assumption \ref{A:1} guarantees that the solution $X(t)$ is a stationary Markov process. This means that the conditional distribution of $X_{t_2}$ with respect to (w.r.t.) its location at a time $t_1<t_2$, only depends on the time elapsed $t_2 - t_1$.
	
	A Markov process can be alternatively studied by the continuous transformation it induces on probability measures. Recall that the distribution of $X_t$ is determined by that of $X_0$. If the distribution of $X_0$ is absolutely continuous w.r.t. the Lebesgue volume measure on $\real^d$, then the same will be true for $X_t$ for every $t>0$. Thus an SDE induces a dynamics on the space of probability densities as well. This dynamics or continuous transformation of the densities can be tracked by the \emph{forward Kolmogorov operator} $\calL$. This a second order differential operator whose action on a $C^2$ \footnote{Recall that a function is called $C^r$ if it is $r$-times differentiable, and its first $r$ derivatives are continuous.} function $g$ is given by
	\begin{equation} \label{eqn:def:FP}
		\paran{ \calL g } (X,t) := \sum_{i=1}^{d} V(x)_i \frac{\partial}{\partial_i} g(x) + \sum_{i,j=1}^{d} \paran{ \frac{\partial^2}{ 2\partial_i \partial_j } g }(x) A_{i,j}(x) ,
	\end{equation}
	where $A(x) = 0.5 G(x) G(x)^T$ is the $d\times d$ covariance matrix. One of the pillars of the theory of Markov processes is that the derivative of the conditional expectations of a stationary Markov process is provided by the generator $\calL$ (see \cite[Cor 160]{ShaliziKontorovich2010}) :
	\begin{equation} \label{eqn:deriv_exp:1}
		\frac{d}{ds} \rvert_{s=t} \Expect{}{ \phi \paran{ X(s) } \,\rvert\, X(t) } = \calL X(t), \quad \forall \phi \in \dom(\calL) .
	\end{equation}
	Equation \eqref{eqn:deriv_exp:1} is a combination of a second order differential operator, the measure theoretic operation of conditional expectation, and a time derivative. It provides a bridge between statistical properties, (i.e. the distribution of the data), and differential properties (i.e., the drift). Of special interest to us is the particular instance of \eqref{eqn:deriv_exp:1} when $\phi = \Id_{\real^d}$ and $t=0$. In that case we get
	\begin{equation} \label{eqn:deriv_exp:3}
		\lim_{t\to 0^+} \frac{1}{t} \Expect{}{ X(t, x_0) - x_0 } = \paran{\calL \Id_{\real^d}}(x_0) = V(x_0) .
	\end{equation}
	Equation \eqref{eqn:deriv_exp:3} is the core mathematical principle of our numerical technique. In fact we use a third order version of \eqref{eqn:deriv_exp:3} given by
	\begin{equation} \label{eqn:drift:Taylor:3}
		\Delta t V(x_0) = 18 \Expect{}{ X \paran{ \Delta t, x_0} - x_0 } - 9 \Expect{}{ X \paran{ 2\Delta t, x_0} - x_0 } + 2 \Expect{}{ X \paran{ 3\Delta t, x_0} - x_0 } + \bigO{ \Delta t^3 } .
	\end{equation}
	This identity is based on a Taylor series expansion of operator $\calL$ :
	\[ \frac{1}{\Delta t} \SqBrack{ \phi(X_{t+\Delta t}) - \phi(X_t) } = \sum_{n} \frac{1}{n!} \calL^n \phi (X_t) \paran{ \Delta t }^n .  \]
	The precise nature of the convergence depends upon the regularity of $V$ and $G$. See \cite[Eqn 21]{Stanton1997nonpar} for more details. The conditional expectations here are taken w.r.t the path-space $\Omega_{\text{path}}$.
	Equations \eqref{eqn:deriv_exp:3} and \eqref{eqn:drift:Taylor:3} provide a connection between successive states of an SDE path, and the drift term of the SDE. Thus the task of estimating the drift becomes the task of estimating certain conditional expectations. This means that the determination of a deterministic, pre-defined vector valued function can be placed in a probabilistic context, and solved by probabilistic means.
	
	A careful design of a procedure for this task requires a careful understanding of the nature of the data. Formally, we assume the following about the data :
	
	\begin{Assumption} \label{A:2}
		There is sampling interval $\Delta t>0$, a collection of states $x_1, \ldots, x_N$ equidistributed with respect to the stationary measure $\mu_{stat}$, and a collection of $N$ 4-tuples of points 
		\begin{equation} \label{eqn:def:sample}
			\SetDef{ \paran{ X \paran{ 3\Delta t, x_n}, X \paran{ 2\Delta t, x_n}, X \paran{ \Delta t, x_n}, x_n } }{ 1\leq n \leq N  }
		\end{equation}
	\end{Assumption}

	Thus our assumption on the data is that we have snapshots of sample trajectories of length at least $4$, and that the data is equidistributed w.r.t. the unknown stationary measure. The assumption of \emph{equidistribution} is essential and fundamental to most data-driven methods for discovering dynamical systems. Our assumption on data requires neither a single orbit, nor long sequences of data. It does allow the points $x_n$ to be gathered from a single trajectory of the SDE, observed at various time instants.  This makes our technique applicable to systems which are being intermittently observed.

	Our goal can now be precisely stated to be the estimation of the drift $V$ from \eqref{eqn:sde}, while working under the constraints and conditions of Assumptions \ref{A:1} and \ref{A:2}. While there has been many a wide range of techniques for approximating the drift for ODEs \cite[e.g.]{gouesbet1991reconstruction, tsutsumi2022constructing, gouesbet1994global, cao2011robust, ait2010operator}, it is still a challenge to develop a technique that fulfills the following criterion :
	\begin{enumerate} [(i)]
		\item The technique is nonparametric; i.e. it does not assume a prescribed parametric format, such as Gaussian processes \cite[e.g.]{Garcia2017nonpar} or semi-parametric forms such as linear diffusion terms \cite[e.g.]{Devlin2019opt}
		\item The technique is convergent with data.
		\item The technique does not assume a prior distribution; such as Gaussian priors \cite[e.g.]{darcy2023one, Batz2018approx}.
		\item The technique is adaptable to high dimensional systems.
	\end{enumerate}
	Effective nonparametric techniques have been built \cite{RuttorEtAl2013approximate, ella2024nonparametric} which are well suited to sparse sampling but rely on the assumption of prior distributions. Another important work \cite{davis2022est} considers data collected from SDE sampled at uneven intervals. Besides the probabilistic approach presented in this article, four other notable and distinct techniques deserve special mention -- a model reduction technique of Ye et. al. \cite{ye2024nonlinear}; the one based on the \emph{Magnus expansion} \cite{Wang2020magnus}; a back-stepping technique that converts the nonlinear process into an OU process \cite{yin2022backstepping}; and the \emph{multi-level parameter estimation network} for parametric SDEs with jump noises. Our technique achieves all of the first three objectives and overcomes the challenge of dimensionality in certain structured vectored fields, as explained later in Section \ref{sec:sparse}. 
	
	\begin{figure}[!t]\center
		\includegraphics[width=0.3\linewidth, height=0.3\textheight, keepaspectratio]{\figs 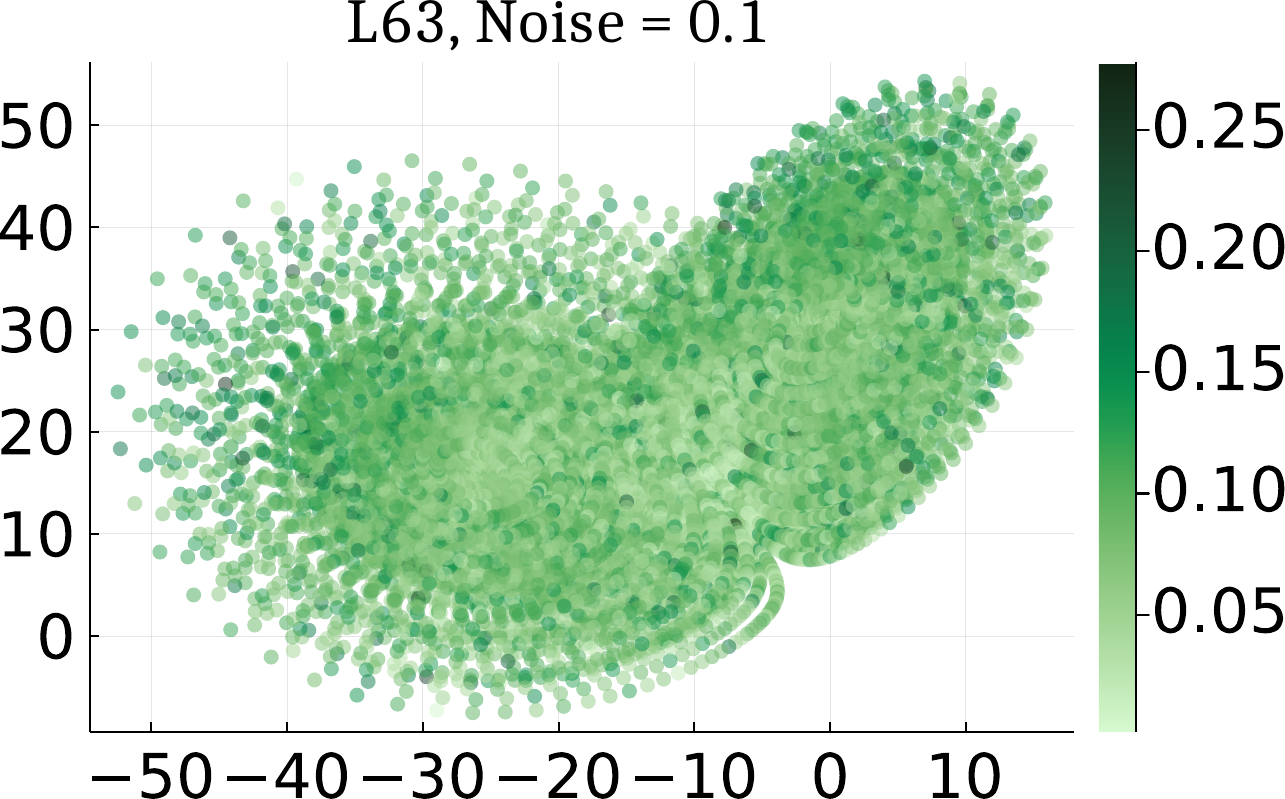}
		\includegraphics[width=0.3\linewidth, height=0.3\textheight, keepaspectratio]{\figs 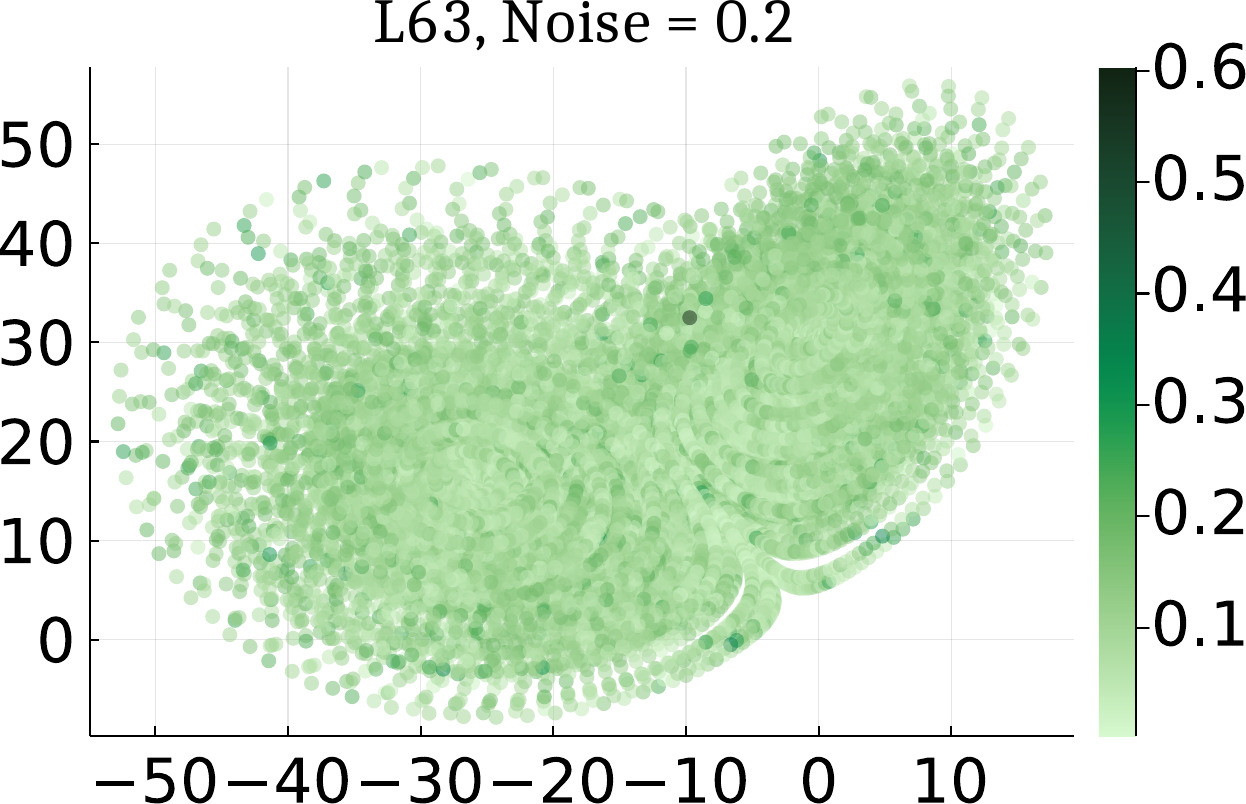}
		\includegraphics[width=0.3\linewidth, height=0.3\textheight, keepaspectratio]{\figs 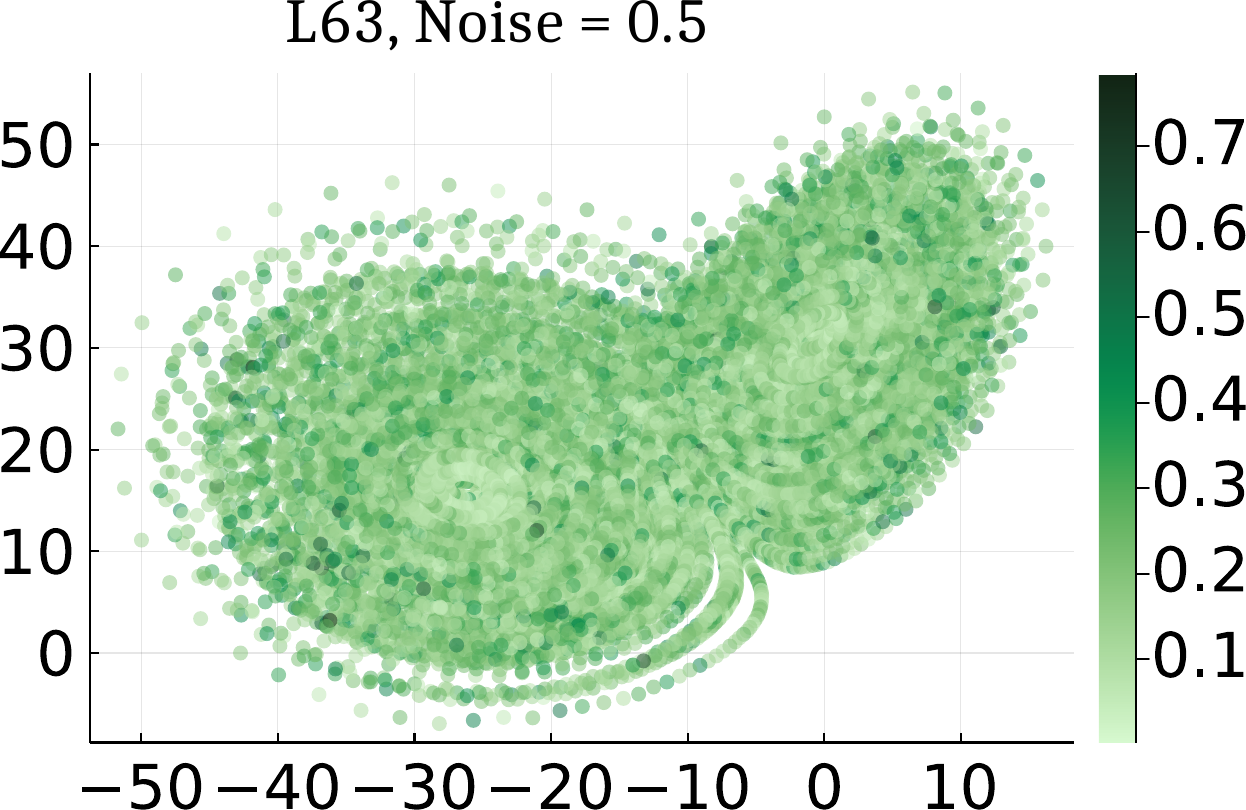}
		\caption{Estimating the drift in a stochastic Lorenz 63 model. The three figures present a heatmap of the pointwise error in computing the drift \eqref{eqn:vec:l63}. The colored area of the phase space represents the support of the stationary measure of the underlying SDE. }
		\label{fig:l63}
	\end{figure}
	
	Our method is based on a measure theoretic insight. Set two spaces 
	\[ \calX = \real^d, \quad \calY = \real^{3d} . \]
	These two spaces $\calX$ and $\calY$ respectively represent the spaces from which $x_0$ and the \\triples $\paran{ X(\Delta t, x_0), X(2\Delta t, x_0), X(3\Delta t, x_0) }$ are drawn. Define a function
	\begin{equation} \label{eqn:def:f}
		f : \calX \times \calY \to \real^d, \quad x_0 \times \paran{ x_1, x_2, x_3 } \mapsto \frac{1}{\Delta t} \SqBrack{ 18 x_1 - 9 x_2 + 2x_3 - 11 x_0} 
	\end{equation}
	Due to the stationarity assumption from Assumption \ref{A:1}, there is a natural measure $\mu$ on $\calX \times \calY$ representing the joint distribution 
	\[ \paran{ x_0, X(\Delta t, x_0), X(2\Delta t, x_0), X(3\Delta t, x_0) } \sim \mu \]
	Let $\proj_1 : \calX \times \calY\to \calX$ denote the projection map to the first component. 
	Then \eqref{eqn:drift:Taylor:3} can be interpreted to say :
	\begin{equation} \label{eqn:deriv_exp:4}
		V = \Expect{\mu}{ f \,\rvert\, \proj_1  } + O \paran{ \Delta t^3 }
	\end{equation}
	Thus the task of estimating the drift has been transformed into the task of finding the conditional estimation of one projection function with respect to another. The conditional expectation is reconstructed or learnt using an operator theoretic tool from \cite{Das2023conditional}. It will be a combination of kernel based techniques and the third order approximation of \eqref{eqn:deriv_exp:4}.  This completes a description of the problems we aim to solve, and the underlying mathematical principles. 
	
	
	\paragraph{Outline} We next present in Section \ref{sec:numerics} a numerical recipe for implementing the conditional expectation route to drift estimation. Next in Section \ref{sec:cnvrg}, we discuss the convergence of our methods, and identify the dependence of the convergence on factors such as data size and noise-levels. Next in Section \ref{sec:examples} we apply the technique to some common examples from Dynamical systems theory. The examples reveal the challenges posed by high dimensional systems. We propose some adjustments to the algorithms to overcome this practical challenge, in Section \ref{sec:sparse}. We end with some discussions in Section \ref{sec:conclus}, summarizing our results and discussing new challenges.
	
	\section{The technique} \label{sec:numerics}
	
	Our technique makes multiple and different uses of kernels and kernel integral methods. A \emph{kernel} on a space $\calZ$ is a bivariate function $k : \calZ \times \calZ \to \real$, and is interpreted to be a measure of similarity between points on $\calZ$. Kernel based methods offer a non-parametric approach to learning, and have been used with success in many diverse fields such as spectral analysis \cite{DasGiannakis_delay_2019, DasGiannakis_RKHS_2018}, 
	discovery of periodic and chaotic components of various real world systems \cite[e.g.]{DasMustAgar2023_qpd, DasEtAl2023traffic}, and even abstract operator valued measures \cite{DGJ_compactV_2018}. The kernels we used are based on \emph{Gaussian} kernels, defined as
	\begin{equation} \label{eqn:def:GaussK}
		k_{\text{Gauss}, \epsilon} (x,y) := \exp \paran{ -\frac{1}{\epsilon} \dist( x, y )^2 } , \quad \forall x,y \in \calZ. 
	\end{equation}
	The bandwidth $\epsilon$ controls how quickly the value of the kernel decays to zero as $x$ and $x'$ grow apart.  We use a Markov normalized version of the Gaussian kernel \eqref{eqn:def:GaussK}. Given a measure $\beta$ on $\calZ$ and a bandwidth parameter $\epsilon>0$, we can create a new kernel
	\begin{equation} \label{eqn:def:GaussSymm}
		k_{ \text{Gauss}, \epsilon }^{ \beta } : \calZ \times \calZ \to \real , \quad k_{ \text{Gauss}, \epsilon }^{ \beta }(x, x') := k_{ \text{Gauss}, \epsilon } (x, x') / \int_{\calZ} k_{ \text{Gauss}, \epsilon } (x, x'') d\beta(x'').
	\end{equation}
	The kernel $k_{ \text{Gauss}, \epsilon }^{ \beta }$ has the property that for every $x$, the integral w.r.t. the second variable $x'$ is $1$. Thus the kernel $k_{ \text{Gauss}, \epsilon }^{ \beta }$ mimics a discrete time Markov process whose state space is $\calZ$. Closely associated to kernels are kernel integral operators (k.i.o.). Now suppose that $\calX$ is a $C^r$-smooth manifold and $k$ is a $C^r$-smooth kernel. Given a probability measure $\alpha$ on $\calZ$, one has the following integral operator associated to $k$ :
	\[ K^\alpha : L^2(\alpha) \to C^r(\calX) , \quad (K^\alpha \phi)(x) := \int_{\calX} k(x,y) \phi(y) d\nu(y) . \]
	Thus if the kernel $k$ is $C^r$, it converts generally non-smooth $L^2(\alpha)$ functions into $C^r$ functions. For this reason, k.i.o.-s are also known as smoothing operators. We also require our kernel to be \emph{strictly positive definite}, which means that given any distinct points $x_1, \ldots, x_N$ in $\calZ$, numbers $a_1, \ldots, a_N$ in $\real$, the sum $\sum_{i=1}^{N} \sum_{i=1}^{N} a_i a_j k(x_i, x_j)$ is non-negative, and zero iff all the $a_i$-s are zero. 
	
	We next present an useful notation : Let $\theta>0$ be a fixed constant. Given a generic linear regression problem $Ax = y$, the $\theta$-regularized least squares solution will be the vector
	\begin{equation} \label{eqn:def:ridge}
		x_{ls, \, \theta-reg} := \SqBrack{ A^T A + \theta }^{-1} A^T y .
	\end{equation}
	We now present our first algorithm from \cite{Das2023conditional} for computing conditional expectation from data. 
	
	\begin{figure}[!t]\center
		\includegraphics[width=0.32\linewidth, height=0.3\textheight, keepaspectratio]{\figs 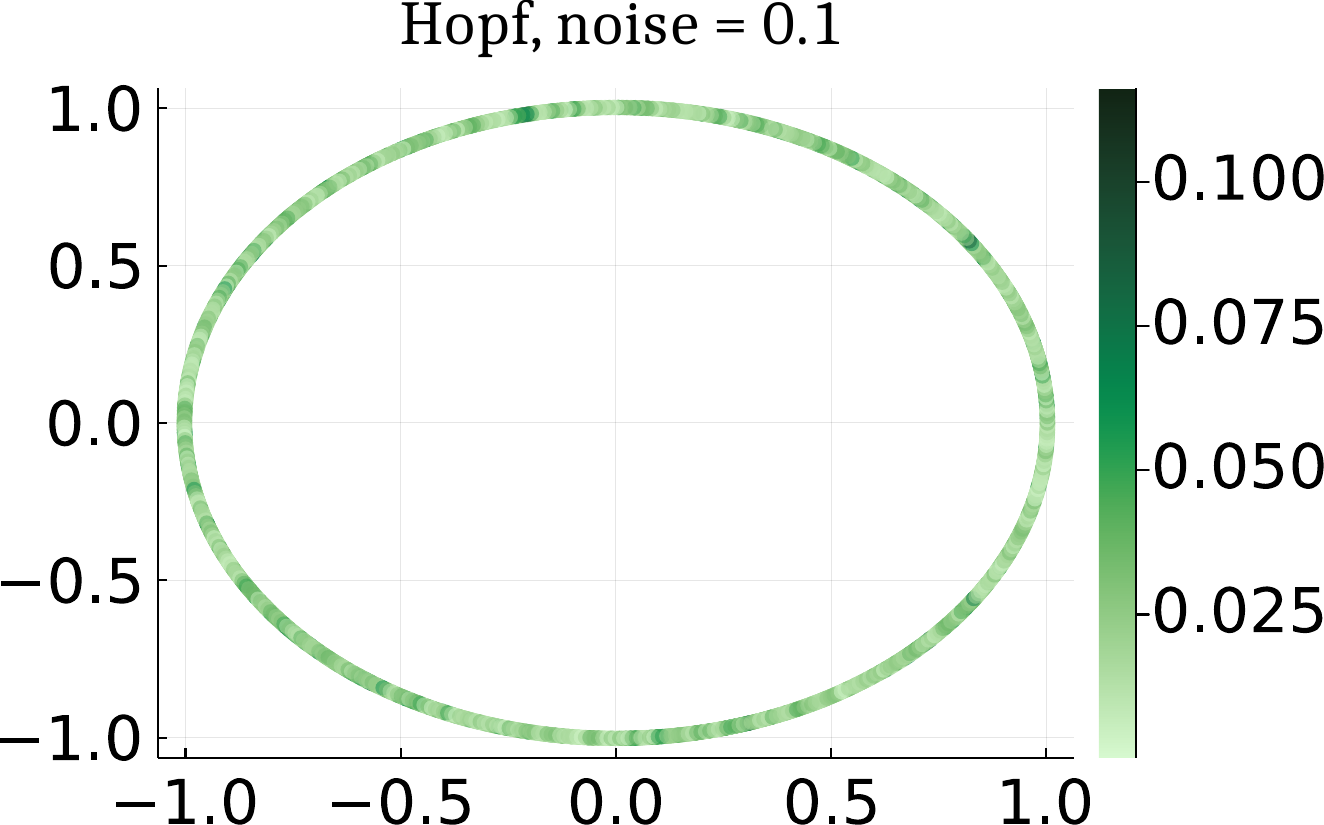}
		\includegraphics[width=0.3\linewidth, height=0.3\textheight, keepaspectratio]{\figs 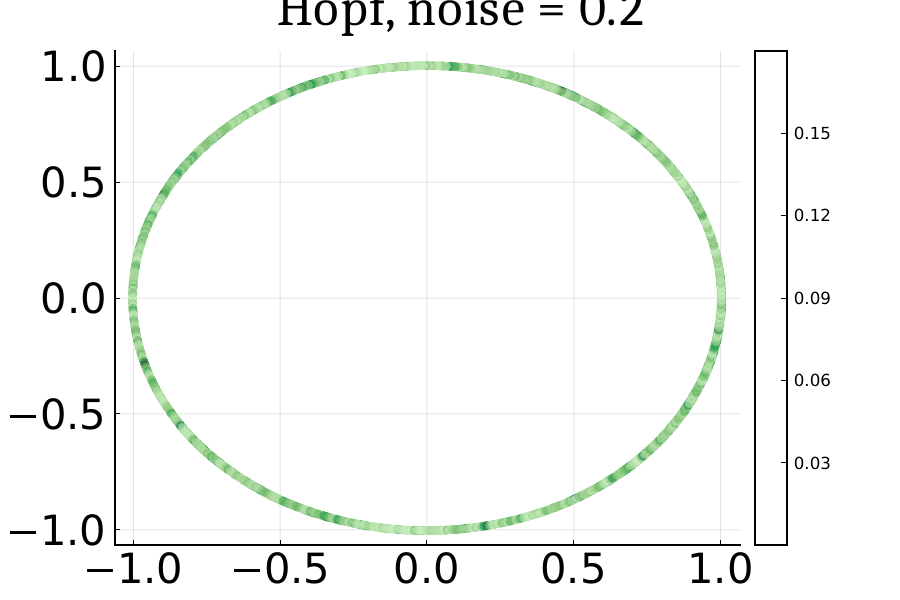}
		\includegraphics[width=0.27\linewidth, height=0.3\textheight, keepaspectratio]{\figs 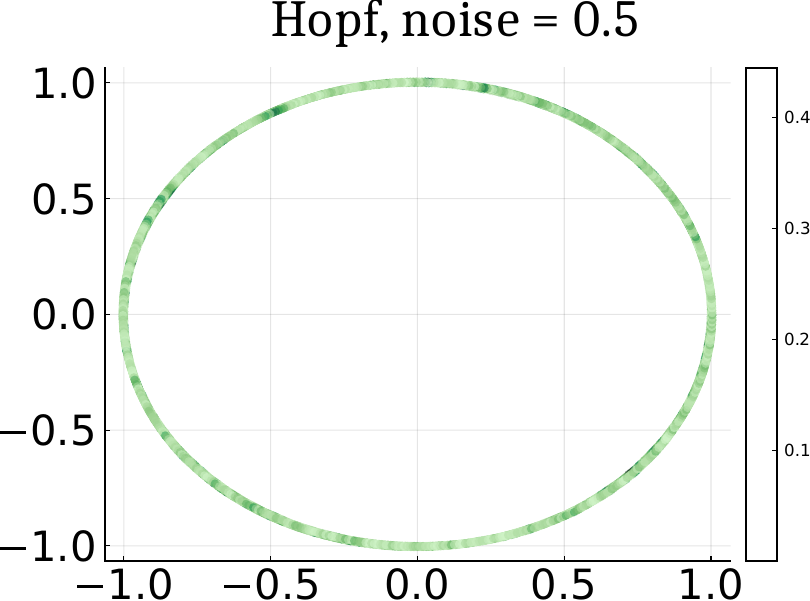}
		\caption{Estimating the drift in a stochastic version of the Hopf oscillator \eqref{eqn:Hopf}. The figures carry the same interpretation as Figure \ref{fig:l63}.}
		\label{fig:Hopf}
	\end{figure}
	
	\begin{algo} \label{algo:1}
		Kernel base estimation of conditional expectation.
		\begin{itemize}
			\item \textbf{Input.} A sequence of pairs $\SetDef{(x_n, y_n)}{ n=1,\ldots,N}$ with $x_n\in\real^d$ and $y_n\in \real$.
			\item \textbf{Parameters.}
			\begin{enumerate}
				\item Strictly positive definite kernel $k:\real^d \times \real^d \to \real^+$.
				\item Smoothing parameter $\epsilon_1>0$.
				\item Sub-sampling parameter $M\in\num$ with $M<N$.
				\item Ridge regression parameter $\theta$.
			\end{enumerate}
			\item \textbf{Output.} A vector $\vec{a} = \paran{a_1, \ldots , a_M} \in \real^M$ that creates a function
			\[ \phi(x) \approx \sum_{m=1}^{N} a_m k(x, x_m) , \quad \forall x \in \real^d , \]
			which approximates the conditional expectation of the y-variable with respect to the x variable.
			\item \textbf{Steps.}
			\begin{enumerate}
				\item Compute a Gaussian Markov kernel matrix using \eqref{eqn:def:GaussSymm} :
				\begin{equation} \label{eqn:GM}
					\Matrix{G_{\epsilon_1}} \in \real^{N\times N}, \quad \Matrix{G_{\epsilon_1}}_{i,j} = k_{ \text{Gauss}, \epsilon }^{ \beta } (x_i, x_j) .
				\end{equation}
				\item Compute a Markov kernel $\Matrix{P} \in \real^{N\times N}$ as $\Matrix{P}_{i,j} := p(x_i, x_j)$
				\item Compute the kernel matrix $\Matrix{K}\in \real^{N\times M}$ as  $\Matrix{K}_{i,j} = k(x_i, x_j)$.
				\item Find a vector $\vec{a} \in \real^M$ as the $\theta$-regularized least-squares solution to the equation
				\begin{equation} \label{eqn:scheme:3}
					\Matrix{P} \Matrix{K} \vec{a} = \Matrix{P} \Matrix{G_{\epsilon_1}} \vec{y} .
				\end{equation}
			\end{enumerate}
		\end{itemize}
	\end{algo}
	
	Algorithm~\ref{algo:1} has two components, the choice of a reconstruction kernel $k$, and a Markov kernel $p$ which approximates the smoothing operator. We usually choose $p$ to be the Markov normalized Gaussian kernel from \eqref{eqn:def:GaussSymm}.  The algorithm is convergent in the following sense : 
	
	\begin{lemma} \label{lem:dj0k3}
		(\cite[Thm 2]{Das2023conditional})  Suppose there is a probability space $(\Omega, \mu)$ and two measurable functions $X : \Omega \to \real^d$ and $Y : \Omega \to \real$. Suppose that the sequence $\braces{ (x_n, y_n) }_{n\in\num}$ is equidistributed with respect to some measure on $\real^d\times \real$. Let the conditional expectation $V = \mathbb{E}_{\mu} (Y|X)$ lie within the reproducing kernel Hilbert space spanned by $k$.  Then the output $\vec{a}$ of Algorithm \ref{algo:1} satisfies :
		\begin{equation} \label{eqn:thm2}
			\lim_{M\to \infty} \lim_{N\to\infty, \theta \to 0^+} \norm{ \sum_{n=1}^{N} a_n k(\cdot, x_n) - \Smooth_{\epsilon_1}^\mu V }_{\rkhs} = 0, \quad \lim_{\epsilon_1 \to 0+} \norm{ \Smooth_{\epsilon_1}^\mu V - V }_{L^2(\nu)} = 0.
		\end{equation}
	\end{lemma}
	
	\begin{figure}[!t]\center
		\includegraphics[width=0.3\linewidth, height=0.3\textheight, keepaspectratio]{\figs 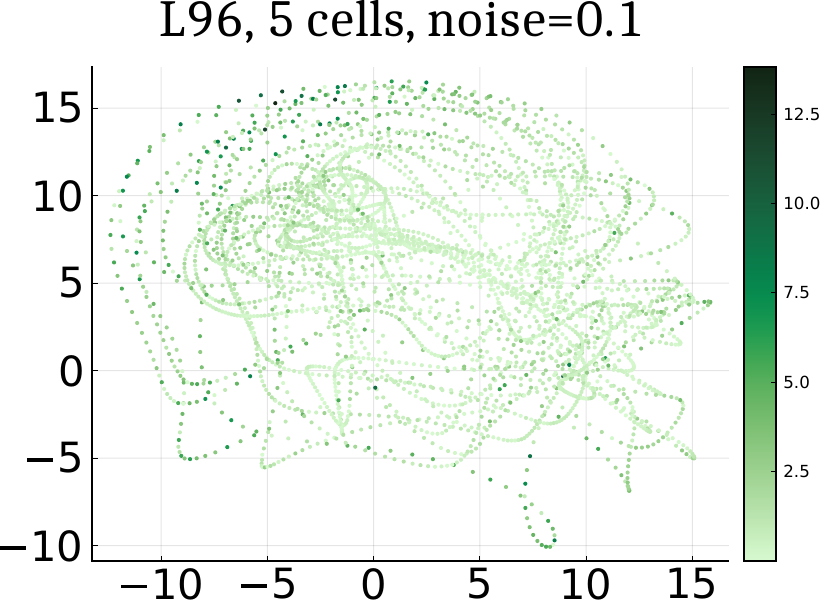}
		\includegraphics[width=0.3\linewidth, height=0.3\textheight, keepaspectratio]{\figs 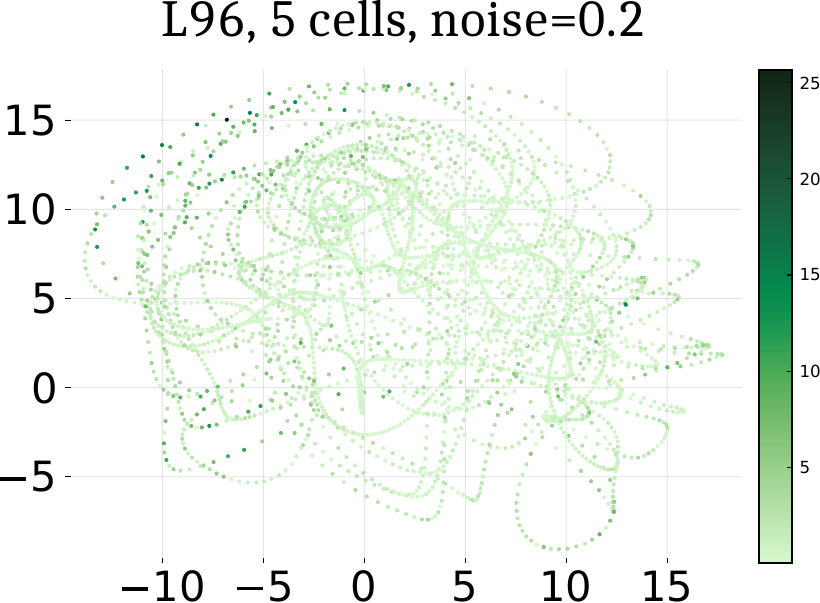}
		\includegraphics[width=0.3\linewidth, height=0.3\textheight, keepaspectratio]{\figs 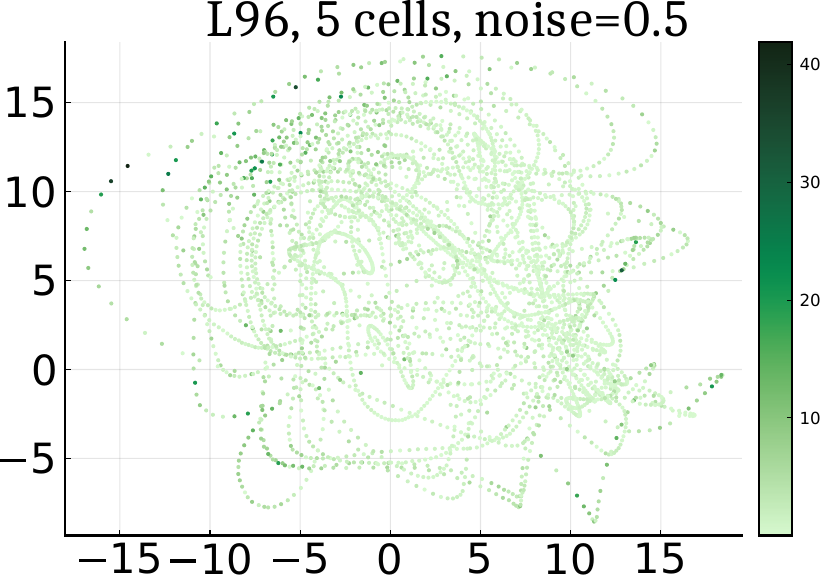}
		\caption{Estimating the drift in a stochastic version of the L96 system \eqref{eqn:L96} with 5 oscillators. The figures carry the same interpretation as Figure \ref{fig:l63}.}
		\label{fig:L96_N5}
	\end{figure}
	
	Algorithm~\ref{algo:1} is a general algorithm for computing conditional expectation, and is at the core of our numerical procedure. We make an important note here that the application of the Gaussian smoothing $\Matrix{G_{\epsilon_1}}$ is optional, and does not affect the analytic properties of the convergence scheme.  We next show its use in drift estimation. Given any $d\times M$ matrix $A$, we use $A_i$ and $A^j$ to denote its $i$-th column vector and $j$-th row vector respectively. Given a sequence $\braces{ z_n }_{n=1}^{N}$ of $d$-dimensional vectors, we use $\braces{ z_n^{(i)} }_{n=1}^{N}$ to denote the $1$-dimensional timeseries created by extracting the $i$-th coordinate from these vectors. We now present how Algorithm~\ref{algo:1} can be used to estimate the drift from an SDE orbit.
	
	\begin{algo} \label{algo:2}
		Drift estimation as a conditional expectation
		\begin{itemize}
			\item \textbf{Input.} Dataset $\braces{ \paran{ X \paran{ 3\Delta t, x_n}, X \paran{ 2\Delta t, x_n}, X \paran{ \Delta t, x_n}, x_n } }_{n=1}^{N}$ as in Assumption \ref{A:2}.
			\item \textbf{Parameters.} Same as Algorithm \ref{algo:1}.
			\item \textbf{Output.} A matrix $A\in \real^{d\times M}$ representing a drift 
			\[ \hat{V}(X) \approx \sum_{m=1}^{M} A_m k(x, x_m) , \quad \forall x\in \real^d. \]
			\item \textbf{Steps.} 
			\begin{enumerate}
				\item Compute 
				\[ z_n = \frac{1}{\Delta t} \SqBrack{ 18 X \paran{ \Delta t, x_n} - 9 X \paran{ 2\Delta t, x_n} + 2 X \paran{ 3\Delta t, x_n} - 11 x_{n} } , \quad \forall n\in 1,\ldots, N . \]
				\item For each $i\in 1,\ldots, d$, pass the timeseries $\braces{ \paran{x_n, z_n^{(i)}} }_{n=1}^{N}$ into Algorithm \ref{algo:1}, along with the rest of the parameters. 
				\item Store the result in a $m$-dimensional column vector $\vec{a}^{(i)}$.
				\item Compute the $d\times M$ matrix $A$ such that
				\[ A^T = \SqBrack{ \vec{a}^{(1)}, \vec{a}^{(2)}, \ldots, \vec{a}^{(d)} } . \]
			\end{enumerate}
		\end{itemize}
	\end{algo}
	
	The convergence of Algorithm \ref{algo:2} is guaranteed by Lemma \ref{lem:dj0k3} and \eqref{eqn:deriv_exp:4}. We investigate this convergence more closely in the next section.
	
	\begin{figure}[!t]\center
		\includegraphics[width=0.95\linewidth, height=0.3\textheight, keepaspectratio]{\figs 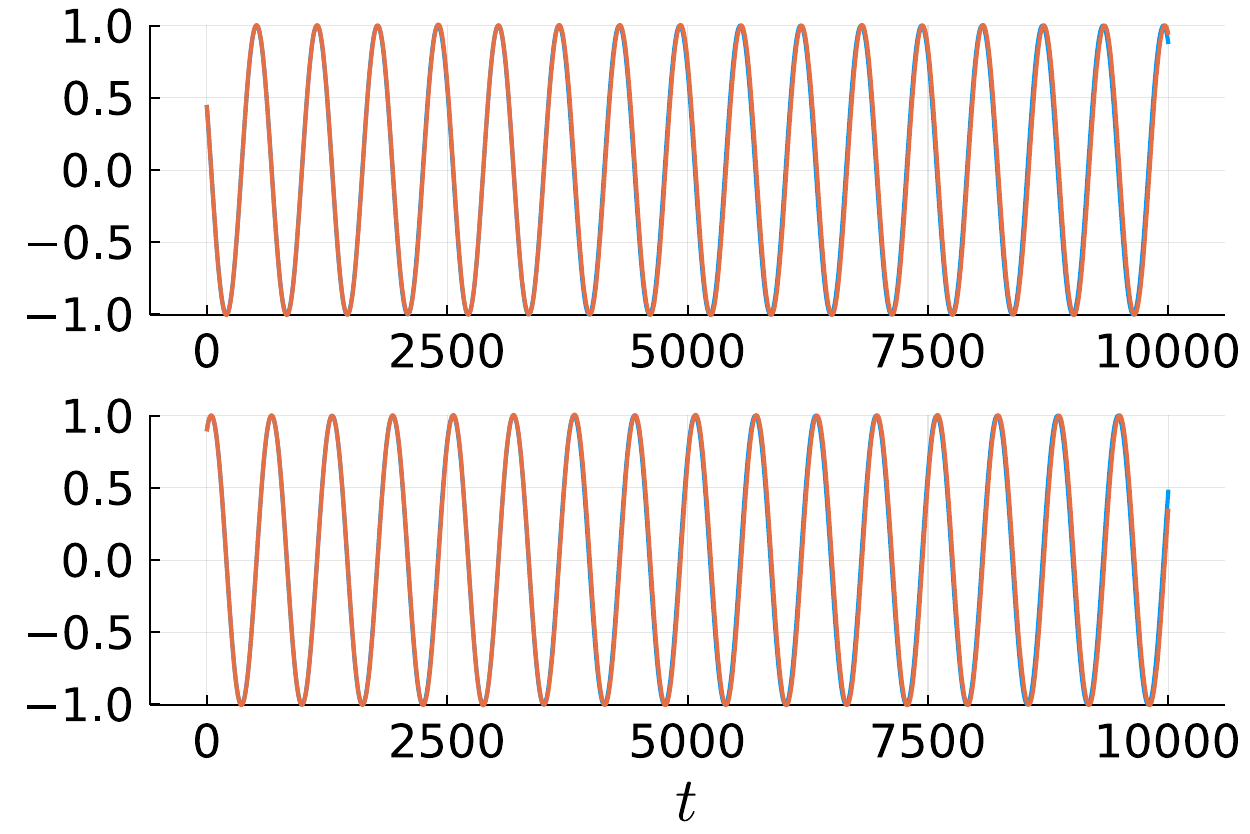}
		\caption{Reconstruction of the Hopf oscillator. The trajectories corresponding to the true and reconstructed drifts are in blue and orange respectively. Due to a strong match in the reconstructed drifts, and the stability of rotational dynamics, the orbits overlap.}
		\label{fig:Hopf_recon}
	\end{figure}
	
	\section{Convergence analysis} \label{sec:cnvrg} 
	
	In this section we conduct a precise analysis on the convergence of our results, and its dependence on various algorithmic and data-source parameters. Throughout this section we assume that Assumptions \ref{A:1} and \ref{A:2} hold; the drift and diffusion terms in \eqref{eqn:sde} are continuous functions; and that $a_{N, M, \theta}$ is the output of Algorithm~\ref{algo:2} with $N$ data points, sub-sampling parameter $M$, and regression regularization parameter $\theta$. We use $\hat{V}_{N, M, \theta}$ to denote the resulting approximation of the drift.
	
	\begin{theorem} [Almost sure convergence] \label{thm:1}
		Suppose Assumptions \ref{A:1} and \ref{A:2} hold. If the drift $V$ lies in the RKHS $\rkhs$ then the following convergence holds almost surely :
		\begin{equation} \label{eqn:scheme:7}
			\boxed{ \lim_{M\to \infty} \lim_{N\to \infty, \theta\to 0^+} \norm{ \hat{V}_{N, M, \theta} - V }_{\sup} \stackrel{a.s.}{=} 0. }
		\end{equation}
	\end{theorem}
	
	Theorem \ref{thm:1} will be shown to be a consequence of functional analytic results from \cite{Das2023conditional}. The convergence result does not reveal the dependence of the parameters to get arbitrarily close to convergence. It also relies on the assumption that $V$ lies in the hypothesis space $\rkhs$. The three main parameters involved are  - 
	\begin{enumerate} [(i)]
		\item the number of data-samples $N$, which dictates how well the distribution of the data points in $\real^{4d}$ is approximated;
		\item choice of the sub-sampling parameter $M$;
		\item a regularization parameter $\theta>0$, as explained in  \eqref{eqn:def:ridge}, for the linear regression problem.
	\end{enumerate} The next theorem presents a more nuanced view of the convergence, in terms of error bounds :

	\begin{theorem} [Parameter selection within error bounds] \label{thm:2}
		Suppose Assumptions \ref{A:1} and \ref{A:2} hold. Fix error bounds $\delta_1, \delta_2, \delta_3 >0$ and an uncertainty bound $\delta_4>0$. 
		\begin{enumerate} [(i)]
			\item For $M$ large enough, $V$ may be approximated uniformly within an error of $\delta_1$ as a sum of kernel sections 
			\[ \norm{ V - \sum_{m=1}^{M} \bar{a}_m k \paran{ \tilde{x}_m, \cdot } }_{ \sup \calX } < \delta_1 . \]
			Here $\braces{ \tilde{x} }_{m=1}^{M}$ is a collection of $M$ points from the dataset assumed in Assumption \ref{A:2}.
			\item Given such a choice of $M$, there exists $\theta_0>0$ and $N'_0\in\num$ such that for every $N>N'_0$ and $\theta \in (0, \theta_0)$, if the diffusion term $G\equiv 0$ , then the approximation error of the approximate drift $\sum_{m=1}^{M} \bar{a}_m k \paran{ \tilde{x}_m, \cdot }$ is less than $\delta_2$.
			\item Fix such a $\theta \in (0, \theta_0)$. Then there is an $N_0$ satisfying
			\[ N_0 > N'_0, \quad N_0 = O \paran{ \frac{d^2}{\delta_3^4 \delta_4^2} } , \]
			such that for every $N>N_0$, with probability at least $1-\delta_4$, the total error in approximation is uniformly bounded by $\delta_1 + \delta_2 + \delta_3$.
		\end{enumerate}\[  ,\]
	\end{theorem}
	
	Theorem \ref{thm:2} also presents the order in which the parameters are to be tuned, starting with $M$ and ending in $N$. The error term $\delta_1$ is related to the approximation error in functional space. The error term $\delta_2$ is controlled by the regression parameter $\theta$. The error term $\delta_3$ is a result of the fluctuations of cumulative error terms from the mean value, as will be described later. The error term $\delta_4$ is just a bound on uncertainty. The gross error $\delta_1 + \delta_2 + \delta_3$ thus reflects an accumulation of errors from three separate levels of approximation. Theorem \ref{thm:2} delineates the contribution of each parameter choice to the error. Note the appearance of $d^2$ in the last lower bound. This term reveals the inevitable challenge of dimensionality. 
	
	The proofs of both the theorems emerge from a careful comparison of the output of the algorithm against the true drift.  Recall the spaces $\calX, \calY$, measure $\mu$ and functions $\proj_1, f$ from Section \ref{sec:intro}. Note that the function $f$ is just the transform on the datapoints, used in Algorithm \ref{algo:2}. Let $\beta$ denote the empirical measure borne by the $N$ data points used in algorithm \ref{algo:2}, and $\beta_X$ denote the projection of $\beta$ to $\calX$. Let $\nu$ denote the empirical measure created from a sub-sampling of $\beta_X$. These notations allow us to interpret the various objects in the numerical procedure in a rigorous operator theoretic setting. With this in mind the regression \eqref{eqn:scheme:3} becomes
	\begin{equation} \label{eqn:scheme:4}
		V \approx \hat{V}_{\beta, \nu, \theta}, \quad \hat{V}_{\beta, \nu, \theta} := K^\nu a_{\beta, \nu, \theta} , \quad a_{\beta, \nu, \theta} := \begin{array}{c} {} \\ \arg\min \\ \theta-\text{reg} \end{array} 
		\SetDef{ \norm{ \mathbb{P}^{\beta_X} K^{\nu} a - \tilde{ \mathbb{P} }^{\beta} f }_{L^2(\nu)} }{ a \in L^2(\nu) } .
	\end{equation}
	Algorithm \ref{algo:1} implements the least-squares regression problem posed in \eqref{eqn:scheme:4}. The measures $\beta$ and $\beta_X$ are completely characterized by the data count $N$, $\nu$ is characterized by a sub-sampling parameter $M$. Recall that given any measure $\alpha$ on $\calX$ we us $P^\alpha$ and $K^\alpha$ to denote the resulting kernel integral operators with kernels $p$ and $k$ respectively. Moreover, we shall use $\tilde{p}$ to denote the trivially extended kernel function
	\[ \tilde{p}: \calX \times \calX \times \calY \to \real, \quad (x_1, x_2, y) \mapsto p(x_1, x_2). \]
	Similarly, given any measure $\tilde{\alpha}$ on $\calX \times \calY$, the resulting integral operator with kernel $\tilde{p}$ will be denoted as $\tilde{P}^{\alpha}$. This notation allows the output $\hat{V}_{\beta, \nu, \theta}$ of Algorithm \ref{algo:1} to be expressed concisely as  
	\begin{equation} \label{eqn:scheme:5}
		\begin{split}
			\hat{V}_{N, M, \theta} := \hat{V}_{\beta, \nu, \theta} &:= \Shobuj{ \Matrix{ K^{\nu} } }
			\akashi{ \SqBrack{ \Matrix{ K^{\nu} }^T \itranga{ \Matrix{ P^{\beta_X}_{\epsilon_3} }^T \Matrix{ P^{\beta_X}_{\epsilon_3} } } \Matrix{ K^{\nu} } + \theta }^{-1} }
			\Shobuj{ \Matrix{ K^{\nu} }^T } \itranga{ \Matrix{ P^{\beta_X}_{\epsilon_3} }^T
				\Matrix{ \tilde{P}_{\epsilon_3}^{\beta} } }
			f \\
			&= \Shobuj{ \Matrix{ K^{\nu} } } \akashi{ \Matrix{ A_{N, \epsilon_3, M} + \theta }^{-1} } \Shobuj{ \Matrix{ K^{\nu} } } \itranga{ \Matrix{ B_{N, \epsilon_3} } } f
		\end{split}
	\end{equation}
	where we have set
	\[ B_{N, \epsilon_3} := \Matrix{ P^{\beta_X}_{\epsilon_3} }^T \Matrix{ P^{\beta_X}_{\epsilon_3} } , \quad A_{N, \epsilon_3, M} := \Matrix{ K^{\nu} } B_{N, \epsilon_3} \Matrix{ K^{\nu} } . \]
	Each of the matrices represented in square brackets above, are also kernel integral operators w.r.t. the finite empirical measures $\beta$, $\nu$ and $\beta_X$.
	Equation \eqref{eqn:scheme:5} reveals that the overall transformation is a sequence of linear transformations. Lemma \ref{lem:dj0k3} and \eqref{eqn:deriv_exp:4} is a direct statement about the almost sure convergence of the $\hat{V}_{\beta, \nu, \theta} = \hat{V}_{N, M, \theta}$ above, thus ensuring the almost-sure convergence \eqref{eqn:scheme:7} expressed in Theorem \ref{thm:1}.
	We now move on to Theorem \ref{thm:2}.
	
	\paragraph{Separating noise} Since we have replaced $\mu$ with $\beta$ and $\mu_X$ with $\nu$, the equality \eqref{eqn:deriv_exp:4} is no longer satisfied. In that case the following difference becomes relevant :
	\begin{equation} \label{eqn:def:err:1}
		\noise_{\beta} := \paran{ \tilde{P}^{\beta} - \tilde{P}^{\mu} }f + \paran{ P^{\beta_X} - P^{\mu_X} } V .
	\end{equation}
	Of course in the limit of infinite data $\beta = \mu$ and $\beta_X = \mu_X$, so the difference $\noise_\beta$ becomes zero. We can now write
	\begin{equation} \label{eqn:def:err:2}
		\tilde{P}^{\beta} f = P^{\beta_X} V + \noise_{\beta}.
	\end{equation}
	The LHS of \eqref{eqn:def:err:2} is the RHS of the regression problem being solved in \eqref{eqn:scheme:4} and \eqref{eqn:scheme:5}. Thus \eqref{eqn:def:err:2} expresses the RHS of the regression problem as an integral transform $P^{\beta_X} V$ of the true field $V$, along with the noise effect $\noise_{\beta}$ defined in \eqref{eqn:def:err:1}. Now note the exact expression for $\noise_{\beta}$ :
	\begin{equation} \label{eqn:def:err:7}
		\paran{\noise_{\beta}}_j = \frac{1}{N} \sum_{n=1}^{N} p(x_j, x_n) \gamma_n, \quad 1\leq j \leq N.,
	\end{equation}
	Recall that each datapoint can be represented as a pair $(x_n, y_n) \in \calX \times \calY$, with $y_n$ having a distribution determined by $x_n$. The sequence $\gamma_n$ are samples of a random variable defined as
	\begin{equation} \label{eqn:def:err:8}
		\gamma_n := f(y_n) - V(x_n).
	\end{equation}
	Equation \eqref{eqn:deriv_exp:4} says that the expected value of the $\gamma$ variable is $O\paran{ \Delta t^3 }$. Equations \eqref{eqn:def:err:1}, \eqref{eqn:def:err:2} and \eqref{eqn:def:err:7} are all equivalent definitions of the noise-residual vector.

	\paragraph{Approximation 1.} Henceforth we shall assume a fixed number $N$ of datapoints $\SetDef{x_n}{1\leq n \leq N}$ distributed evenly over the domain $\calX$, and a subsample $\SetDef{\tilde{x}_n}{1\leq n \leq N}$ of size $M$. This leads to an RKHS subspace
	\[ \rkhs_M := \text{span} \SetDef{ k(\tilde{x}_n, \cdot) }{ 1\leq n \leq M } . \]
	Fix an error level $\delta_1$. Let $B(\calX)$ be the collection of bounded functions on $\calX$. Consider the set
	\[ \rkhs_M(\delta_1) := \SetDef{ \hat{V} \in B(\calX) }{ \exists h\in \rkhs_M \mbox{ s.t. } \norm{ \hat{V} - h }_{B(\calX)} < \delta_1 } . \]
	Thus $\rkhs_M(\delta_1)$ contains those bounded functions which can be uniformly approximated within an error limit of $\delta_1$ by an RKHS function in $\rkhs_M$. Also note that
	\begin{equation} \label{eqn:approx:1}
		\text{clos} \cup_{M\in\num} \rkhs_M(\delta_1) = B(\calX), \quad \forall \delta_1>0.
	\end{equation}
	We shall assume henceforth that the true image $V$ lies in $\rkhs_M(\delta_1)$. Thus there is a vector $\hat{a}\in L^2(\nu)$ such that $\norm{ V - K^\nu \hat{a} } < \delta_1$. This proves Theorem \ref{thm:2}~(i). Allowing an error limit of $\delta_1$ we shall assume without loss of generality that $V = K^\nu \bar{a}$.
	
	\begin{table}
		\caption{Summary of parameters in experiments. In all these experiments, the number of \ldots }
		\begin{tabularx}{\linewidth}{|L|L|L|L|L|}
			\hline
			System & properties & Figures & ODE Parameters & Data samples \\ \hline
			Lorenz 63 \eqref{eqn:vec:l63} & Chaotic system in $\real^3$ & \ref{fig:l63}, \ref{fig:L63_recon} & $\sigma = 10$, $\beta = 8/3$, $\rho=28$ & $10^4$ \\ \hline
			Hopf oscillator \eqref{eqn:Hopf} & Stable periodic cycle in $\real^2$ & \ref{fig:Hopf}, \ref{fig:Hopf_recon} & $\mu=1$ & $10^4$ \\ \hline
			Lorenz 96 \eqref{eqn:L96} & Chaotic system of $10$ 1-dimensional oscillators, cyclically coupled & \ref{fig:L96_N5}, \ref{fig:L96_recon_a} & $f=8.0$, $N=10$ & $2\times 10^3$  \\ \hline
		\end{tabularx}
		\label{tab:param1}
	\end{table}
	
	\paragraph{Approximation 2} Equations \eqref{eqn:scheme:5}, \eqref{eqn:def:err:1} and \eqref{eqn:def:err:2} together imply
	\begin{equation} \label{eqn:def:err:3}
		\begin{split}
			\hat{V}_{N, M, \theta} &= \Shobuj{ K^{\nu} } \akashi{ \SqBrack{ \Matrix{ K^{\nu} }^T \itranga{ \Matrix{ P^{\beta_X}_{\epsilon_3} }^T \Matrix{ P^{\beta_X}_{\epsilon_3} } } \Matrix{ K^{\nu} } + \theta }^{-1} } \Shobuj{ \Matrix{ K^{\nu} }^T } \itranga{ \Matrix{ P^{\beta_X}_{\epsilon_3} }^T } \SqBrack{ P^{\beta_X} V + \noise_{\beta} } \\
			&= \Shobuj{ K^{\nu} } \akashi{ \SqBrack{ A_{N, \epsilon_3, M} + \theta }^{-1} } \Shobuj{ \Matrix{ K^{\nu} } } \SqBrack{ B_{N, \epsilon_3}V + \Matrix{ P^{\beta_X}_{\epsilon_3} }^T \noise_{\beta} } .
		\end{split}
	\end{equation}
	Incorporating the assumption $V = K^\nu \bar{a}$ gives
	\begin{equation} \label{eqn:def:err:4}
		\hat{V}_{N, M, \theta} = \Shobuj{ K^{\nu} } \akashi{ \SqBrack{ A_{N, \epsilon_3, M} + \theta }^{-1} } \akashi{ \SqBrack{ A_{N, \epsilon_3, M} } } \bar{a} + \Shobuj{ K^{\nu} } \akashi{ \SqBrack{ A_{N, \epsilon_3, M} + \theta }^{-1} } \Shobuj{ \Matrix{ K^{\nu} } } \Matrix{ P^{\beta_X}_{\epsilon_3} }^T \noise_{\beta} .
	\end{equation}
	Thus $\hat{V}_{N, M, \theta}$ is the sum of two separate and mutually independent terms, one which depends on the true drift and one which only depends on the noise. We shall examine these terms separately.
	
	\paragraph{Approximation 3} Note that the limit in \eqref{eqn:scheme:7} holds irrespective of the nature of $\mu$. When the diffusion is zero, for each $x\in \calX$ the conditional measures $\mu_{\calY|x}$ are point measures at zero. In that case according to \eqref{eqn:def:err:1} and \eqref{eqn:def:err:8}, $\noise_{\nu} \equiv 0$. Thus \eqref{eqn:scheme:7} and \eqref{eqn:def:err:4} combine to give
	\[\lim_{N\to \infty, \theta\to 0^+} \norm{ \SqBrack{ A_{N, \epsilon_3, M} + \theta }^{-1} \SqBrack{ A_{N, \epsilon_3, M} } \bar{a} - \bar{a} }_{L^2(\nu)} = 0. \]
	Since we keep the sub-sampling parameter $M$ fixed and thus $\nu$ fixed, $K^\nu$ remains a fixed and bounded operator. As a result for fixed $M$ we have :
	\begin{equation} \label{eqn:def:err:5}
		\lim_{N\to \infty, \theta\to 0^+} \norm{ K^{\nu} \SqBrack{ A_{N, \epsilon_3, M} + \theta }^{-1} \SqBrack{ A_{N, \epsilon_3, M} } \bar{a} - K^{\nu} \bar{a} }_{C(\calX)} = 0. 
	\end{equation}
	One important property of the convergence in \eqref{eqn:def:err:5} is that $N, \theta$ can be simultaneously taken to their limiting values. This joint convergence is one of the distinguishing feature of our technique. The limit in \eqref{eqn:def:err:5} proves Theorem \ref{thm:2}~(iii). It remains to prove Theorem \ref{thm:2}~(ii). For that it suffices to bound the error in the second term of \eqref{eqn:def:err:4} by $\delta_3$. In the remainder of the proof our focus will be on this term $\Shobuj{ K^{\nu} } \akashi{ \SqBrack{ A_{N, \epsilon_3, M} + \theta }^{-1} } \Shobuj{ \Matrix{ K^{\nu} } } \Matrix{ P^{\beta_X}_{\epsilon_3} }^T \noise_{\beta}$.
	
	\paragraph{Approximation 4} The self adjoint operator $B_{\beta_{X}, \epsilon_3} = \paran{P^{\beta_{X}}_{\epsilon_3}}^T P^{\beta_{X}}_{\epsilon_3}$ is also a kernel integral operator with kernel
	\[ b^{\beta_{X}}_{\epsilon_3} : \calX \times \calX \to \real, \quad (z, x') \mapsto \int_{\calX} p(z, x') p(x, x') d\beta_{X}(x) \]
	Then according to \eqref{eqn:def:err:2}, we have
	\[\begin{split}
		\SqBrack{ \paran{ P^{\beta_X}_{\epsilon_3} }^T \noise_\beta }(z) &= \int_{\calX} b^{\beta_{X}}_{\epsilon_3}(z, x') \SqBrack{ \int_{\calY} f(x', y) d \beta_{\calY|x'}(y) - V(x') } d\beta_X(x') \\
		&= \int_{\calX} b^{\beta_{X}}_{\epsilon_3}(z, x') \SqBrack{ \int_{\calY} f(x', y) d \beta_{\calY|x'}(y) - \int_{\calY} f(x', y) d \mu_{\calY|x'}(y) } d\beta_X(x')\\
		&= \int_{\calX} b^{\beta_{X}}_{\epsilon_3}(z, x') \int_{\calY} \SqBrack{ f(x', y)- V(x') } d \beta_{\calY|x'}(y) d\beta_X(x') .
	\end{split}\]
	Equation \eqref{eqn:def:err:8} interprets the noise variable $\gamma$ to be additive. Thus we have
	\[ \SqBrack{ \paran{ P^{\beta_X}_{\epsilon_3} }^T \noise_\beta }(z) = \SqBrack{ \paran{ P^{\beta_X}_{\epsilon_3} } \noise_\beta }(z) = \frac{1}{N} \sum_{n=1}^{N} b^{\beta_{X}}_{\epsilon_3}(z, x_n) \gamma_n \]
	The distribution of increments of a process defined by an SDE usually lacks an analytic formula. However due to the continuity of $V$ and $G$ in \eqref{eqn:sde} and the continuity of sample paths of an SDE, we can say that
	\[ \lim_{\Delta t\to 0^+} \text{var}(\gamma) = \mathbb{E}_{\mu_{stat}} GG^T . \]
	In other words, for positive $\Delta t$, $\gamma$ converges in distribution to a Gaussian variable on $\real^d$ whose variance matrix is a strictly positive definite matrix. Now by the central limit theorem we have
	\[\begin{tikzcd}
		X(N,z, \epsilon_3) := N^{1/2} \SqBrack{ \paran{ P^{\beta_X}_{\epsilon_3} } \noise_\beta }(z) \arrow{rr}{ a.s. }[swap]{N\to\infty} &&  \calN \paran{ 0, \text{var}(\gamma) } + O\paran{ \Delta t^3 } .
	\end{tikzcd}\]
	%
	We have declared the quantity to be a random variable $X(N,z, \epsilon_3)$ to denote its dependence on $N$, bandwidth $\epsilon_3$ and initial location $z$.  Let us fix an error bound $\delta_3$. Thus by the multi-variate Chebyshev inequality \cite[e.g.]{MarshallOlkin1960Cheby} we have
	\[ \text{Prob} \braces{ \abs{ P^{\beta_X}_{\epsilon_3} \noise_\beta - O\paran{ \Delta t^3 } }  < \delta_3 } > 1 - O\paran{ \frac{d}{N^{1/2} \delta_3^2} } . \]
	So to bound the uncertainty by $\delta_4$ we would need $N = O \paran{ \frac{d^2}{\delta_3^4 \delta_4^2} } $, as claimed. The completes the proof of Theorem \ref{thm:2}. 
	
	This concludes the description of our methods and the underlying theory. We next describe some numerical experiments to test our technique on.
	\section{Examples} \label{sec:examples} 
	
	In this section we put to test our algorithms. The systems that we investigate are as follows : 
	
	\begin{enumerate}
		\item \textbf{Lorenz 63.} This is the benchmark problem for studying chaotic phenomenon in three dimensions, which is the lowest dimension in which chaos is possible.
		\begin{equation} \label{eqn:vec:l63}
			\begin{split}
				\frac{d}{dt} x_1(t) &= \sigma(x_2-x_1) \\
				\frac{d}{dt} x_2(t) &= x_! (\rho-x_3) - x_2\\
				\frac{d}{dt} x_3(t) &= x_1 x_2 - \beta x_3
			\end{split}
		\end{equation}
		\item \textbf{Hopf oscillator.} The Hopf oscillator is used as a simple parametric model to study the phenomenon of Hopf bifurcations. Given a constant $p$, the ODE is
		\begin{equation} \label{eqn:Hopf}
			\begin{split}
				\frac{d}{dt} x_1(t) &= -x_2 + x_1 \paran{ p- \paran{ x_1^2 + x_2^2 } } \\
				\frac{d}{dt} x_2(t) &= x_1 + x_2 \paran{ p- \paran{ x_1^2 + x_2^2 } }
			\end{split}
		\end{equation}
		\item \textbf{Lorenz 96.} The Lorenz 96 model is a dynamical system formulated by Edward Lorenz in 1996
		\begin{equation} \label{eqn:L96}
			\frac{d}{dt} x_n(t) = \paran{ x_{n+1} - x_{n-2} } x_{n-1} - x_n + F , \quad \forall n\in 1, \ldots, N .
		\end{equation}
		This model mimics the time evolution of an unspecified weather measurement collected at $N$ equidistant grid points along a latitude circle of the earth \cite{LorenzEmanuel1998opti}. The constant $F$ is known as the \emph{forcing} function. The system \eqref{eqn:L96} is benchmark problem in data assimilation.
	\end{enumerate}
	
	\begin{table}[!t]
		\caption{Parameters which remain constant in the experiments. }
		\begin{tabularx}{\linewidth}{|L|L|L|L|} \hline
			Variable & $\eta_3$ & $(l, w)$ & $\theta$ \\ \hline
			Interpretation & Parameter for selecting the bandwidth $\epsilon_3$ of the Markov kernel, based on Algorithm \ref{algo:2} & Dimensions of the image & Ridge regression parameter, as described in \eqref{eqn:def:ridge} \\ \hline
			Value & $4$ & $(1,1)$ & $0.01 \times \norm{K^\nu}$ \\ \hline
		\end{tabularx}
		\label{tab:param2}
	\end{table}
	
	\begin{table}
		\caption{Auto-tuned parameters in experiments. The strategy outlined in \eqref{eqn:gid9c} is used to determine the three bandwith parameters $\epsilon_1, \epsilon_2, \epsilon_3$. }
		\begin{tabularx}{\linewidth}{|L|L|L|L|}
			\hline
			Ridge regression parameter $\epsilon$ & Bandwidth $\epsilon_1$ of Gaussian kernel $\Smooth$ & Bandwidth $\epsilon_2$ of RKHS kernel $k$ & Bandwidth $\epsilon_3$ of Markov kernel $p$ \\ \hline
			0.1 & Chosen to achieve sparsity $0.3\%$ & Chosen to achieve sparsity $1\%$ & Chosen to achieve sparsity $1\%$ \\ \hline
		\end{tabularx}
		\label{tab:param3}
	\end{table}
	
	The steps in each experiments were as follows : 
	\begin{enumerate}
		\item The parameters for the ODEs were set according to Table \ref{tab:param1}.
		\item The ODEs were converted into an SDE with the help of a diagonal diffusion term 
		\begin{equation} \label{eqn:dffsn:expt}
			G(x)_{i,i} = \sigma V(x)_i , \quad i=1,\ldots, d.
		\end{equation}
		\item The constant $\sigma$ which controls the level of noise, was varied among the values $\braces{ 0.1, 0.2, 0.5 }$ respectively. 
		\item Algorithm \ref{algo:2} was applied to a sample path of the SDE. The parameters for the algorithm were auto-tuned as described in Table \ref{tab:param2}.
		\item The choice of kernel $k$ was a diffusion kernel as in \eqref{eqn:def:kdiff}, and a Gaussian Markov normalized kernel as in \eqref{eqn:def:GaussSymm}.
		\item Some of the algorithmic parameters were auto-tuned according to the hyper-parameters listed in Table \ref{tab:param3}.
	\end{enumerate}
	
	The choice of noise levels do not have any specific design, other than providing an increasing sequence of noise levels. Recall that a higher noise level requires a larger number of data points to approximate the conditional expectation. The highest value of $0.5$ was chosen heuristically based on the size $\approx 10^4$ of the dataset being used. We next describe our parameter choices, and methodology for evaluating our results.
	
	\begin{figure}[!t]\center
		\includegraphics[width=0.95\linewidth, height=0.3\textheight, keepaspectratio]{\figs 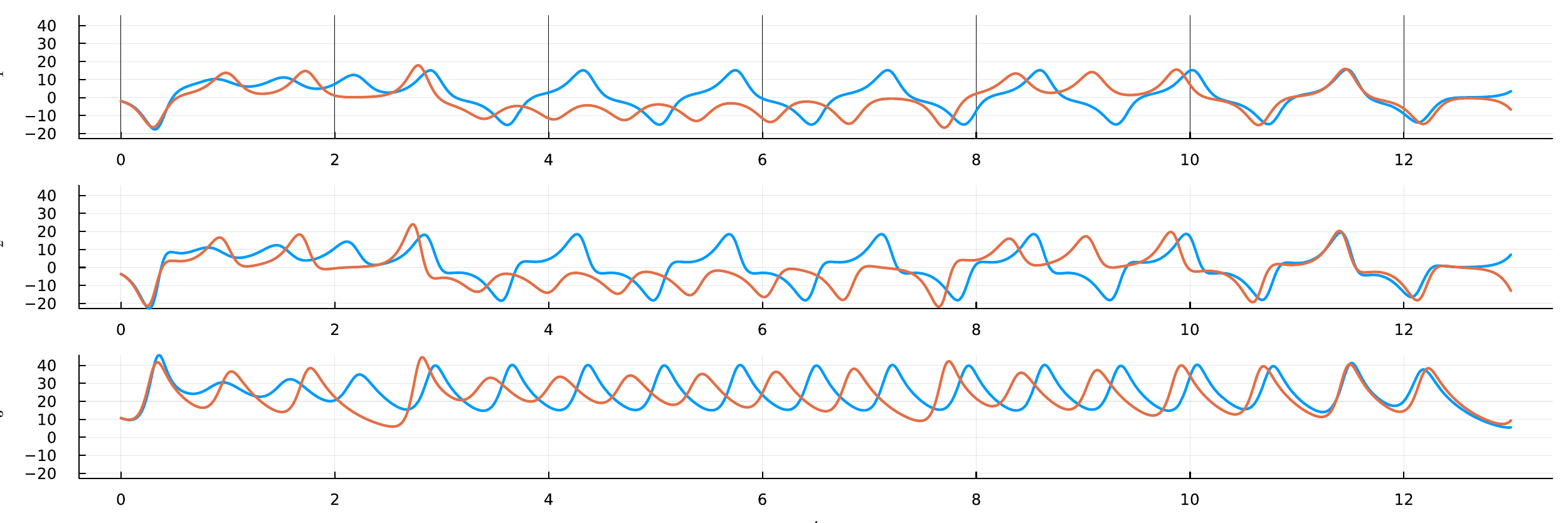}
		\caption{Reconstruction of the Lorenz 63 system. The trajectories corresponding to the true and reconstructed drifts are in blue and orange respectively.}
		\label{fig:L63_recon}
	\end{figure}
	
	\begin{table}
		\caption{Relative $L^2$ error in experiments }
		\begin{tabularx}{\linewidth}{|L|L|L|L|}
			\hline
			Experiment & Noise 0.1 & Noise 0.2 & Noise 0.5  \\ \hline
			Hopf oscillator \eqref{eqn:Hopf} & 0.02 & 0.03 & 0.06 \\ \hline
			Lorenz 63 \eqref{eqn:vec:l63} & 0.08 & 0.1 & 0.14  \\ \hline
			Lorenz 96 \eqref{eqn:L96} - 5 cells & 0.31 & 0.27 & 0.265 \\ \hline
		\end{tabularx}
		\label{tab:results}
	\end{table}
	
	\paragraph{Choice of kernel} Our choice of kernel in all the experiments is the \emph{diffusion} kernel  \cite[e.g.][]{MarshallCoifman2019, WormellReich2021}. There are many versions of the diffusion kernel. We choose the following :
	\begin{equation} \label{eqn:def:kdiff}
		\begin{split}
			& k_{ \text{diff}, \epsilon}^\mu (x,y) = \frac{ k_{\text{Gauss}, \epsilon}(x,y) }{ \deg_l(x) \deg_r(y) } ,\\ 
			& \deg_r(x) := \int_X k_{\text{Gauss}, \epsilon}(x,y) d\mu(y) , \quad \deg_l(x) := \int_X k_{\text{Gauss}, \epsilon}(x,y) \frac{1}{ \deg_r(x) } d\mu(y) .
		\end{split}
	\end{equation}
	Diffusion kernels have been shown to be good approximants of the local geometry in various different situations \cite[e.g.]{CoifmanLafon2006, HeinEtAl2005, VaughnBerryAntil2019}, and are a natural choice for non-parametric learning. We go one step further and perform a symmetrization : 
	\begin{equation} \label{eqn:symm}
		\rho(x) k_{ \text{diff}, \epsilon}^\mu (x,y) \rho(y)^{-1} = \tilde{k}_{ \text{diff}, \epsilon}^\mu (x,y) = \frac{k_{\text{Gauss}, \epsilon}(x,y)}{ \SqBrack{ \deg_r(x) \deg_r(y) \deg_l(x) \deg_l(y) }^{1/2} }, 
	\end{equation}
	where
	\[ \rho(z) = \deg_l(z)^{1/2} / \deg_r(z)^{1/2} . \]
	The kernel $\tilde{k}_{ \text{diff}, \epsilon}^\mu$ from \eqref{eqn:symm} is clearly symmetric. Since it is built from the s.p.d. kernel $k_{ \text{Gauss}, \epsilon}$, $\tilde{k}_{ \text{diff}, \epsilon}^\mu$ is s.p.d. too and thus generates an RKHS of its own. Moreover, the kernel $k_{ \text{diff}, \epsilon}^\mu$  can  be symmetrized by a degree function $\rho$, which is bounded as well as bounded above $0$. Such a kernel will be called \emph{RKHS-like}. Let $M_\rho$ be the multiplication operator with $\rho$. Then
	\[ \ran K^\mu_{ \text{diff}, \epsilon} = \ran M_\rho \circ \tilde{K}^\mu_{ \text{diff}, \epsilon} . \]
	Again, because of the properties of $\rho$, both $M_\rho$ and its inverse are bounded operators. Thus there is a bijection between the RKHS generated by $\tilde{k}_{ \text{diff}, \epsilon}^\mu$, and the range of the integral operator $K_{ \text{diff}, \epsilon}^\mu$. 
	
	\paragraph{Bandwidth selection} The proposed algorithm depends on several choices of bandwidth parameters, which also play the role of a smoothing parameter. We now present our method of selecting a bandwidth $\epsilon$ for constructing a kernel over a dataset $\calD$. It depends on the choice of a threshold $\eta \in (0,1)$, and a zero threshold $\theta_{zero} = 10^{-14}$. We first choose a subsample $\calD'$ of the data set $\calD$ and construct the set $\calS$ of all possible pairwise squared distances. Usually, it suffices to choose $\calD'$ to be 0.1 fraction of the dataset $\calD$, equidistributed throughout $\calD$. Given the threshold $\eta$ , we set $\epsilon$ such that $\eta$-fraction of the numbers in $\calS$ are greater than $\epsilon \theta$, where $\theta$ is the threshold :
	\begin{equation} \label{eqn:gid9c}
		\theta = -1/ \ln \paran{ \theta_{zero} }.
	\end{equation}
	Thus is $\epsilon$ is chosen in this manner an $\eta$ fraction of the set $\SetDef{ e^{ - \norm{x-x'}^2 / \epsilon} }{x,x'\in \calD'}$ are greater than $\theta_{zero}$. The strategy outlined in \eqref{eqn:gid9c} is used to determine the three bandwith parameters $\epsilon_1, \epsilon_2, \epsilon_3$ in Table \ref{tab:param2}.
	
\begin{figure}[!t]\center
		\includegraphics[width=0.95\linewidth, height=0.3\textheight, keepaspectratio]{\figs 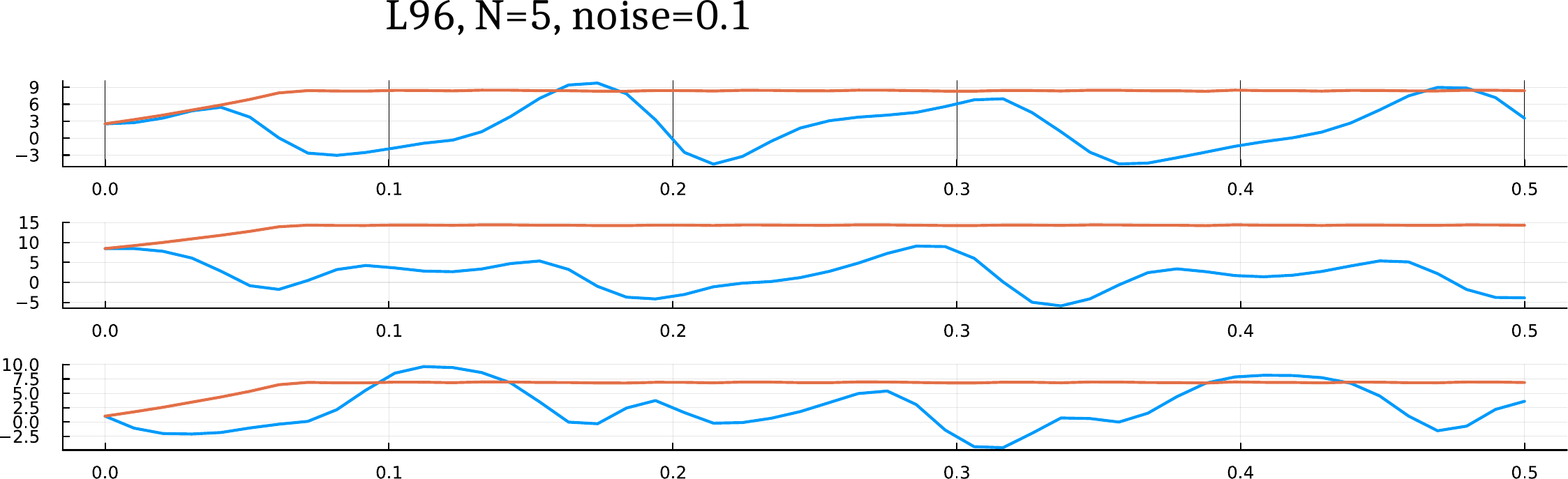}
		\caption{Reconstruction of the Lorenz 96 system with 5 cells. The orbits of the true and reconstructed systems diverge rapidly due to the presence of positive Lyapunov exponents in the L96 system. See Figure \ref{fig:L96_recon_b} for a closer look at the plots. The trajectories corresponding to the true and reconstructed drifts are in blue and orange respectively.}
		\label{fig:L96_recon_a}
	\end{figure}

\begin{figure}[!t]\center
		\includegraphics[width=0.8\linewidth, height=0.3\textheight, keepaspectratio]{\figs 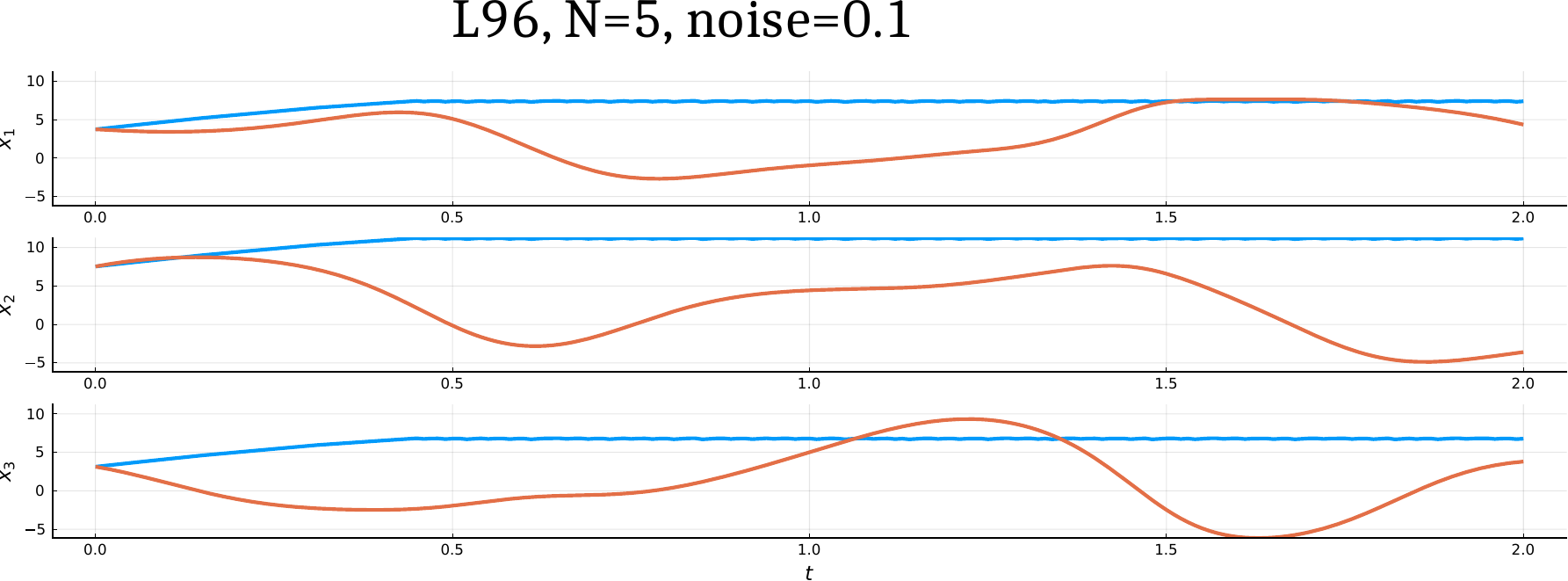}
		\caption{Reconstruction of the Lorenz 96 system with 5 cells. This figure offers a closer look at the graphs of Figure \ref{fig:L96_recon_a}. The trajectories corresponding to the true and reconstructed drifts are in blue and orange respectively.}
		\label{fig:L96_recon_b}
	\end{figure}

	\paragraph{Evaluating the results} There is no certain way of evaluating the efficacy of a reconstruction technique for chaotic systems, as argued in \cite{BerryDas2023learning, BerryDas2024review}. In any data-driven approach such as ours, the reconstruction is done in dimensions higher than the original system. As a result, new directions of instability are introduced, and even if the drift is reconstructed with good precision, the results of their simulation might be completely different. For that purpose, our primary means of evaluating the accuracy of our drift reconstruction is by a pointwise evaluation of the reconstruction error, as shown in Figures \ref{fig:l63}, \ref{fig:Hopf} and \ref{fig:L96_N5}. The average errors are tabulated in  Table \ref{tab:results}. Also see Figures \ref{fig:Hopf_recon}, \ref{fig:L63_recon}, \ref{fig:L96_recon_a} and \ref{fig:L96_recon_b} for a comparison of orbits from the true system and the reconstructed system. It can be noted in the last two figures that the orbits of the reconstructed system get stuck in a false fixed point. This is an artifact of the technique of making out of sample evaluations, which we discuss next.
	
	\paragraph{Evaluating kernel functions} Any kernel function $\phi(x) = \sum_{n=1}^{N} a_n \tilde{K}^\mu_{ \text{diff}, \epsilon} (x, x_n)$ created from a diffusion kernel \eqref{eqn:def:kdiff} can be evaluated at an arbitrary point $x$ in $\real^d$ as 
	\[ \phi(x) = \frac{1}{N} \SqBrack{ \vec{\rho}_l(x) \vec{\rho}_r(x) }^{-0.5} \sum_{n=1}^{N} k_{Gauss, \epsilon} (x,x_n) w_n ,  \]
	where the other terms are given by 
	\[\begin{split}
		w_n = a_n \SqBrack{ \deg_r(x_n) \deg_l(x_n) }^{-0.5} , \; \vec{\rho}_r(x) = \frac{1}{N} \sum_{n=1}^{N} k_{Gauss, \epsilon} (x,x_n) , \; \vec{\rho}_l(x) = \frac{1}{N} \sum_{n=1}^{N} k_{Gauss, \epsilon} (x,x_n) / \deg_r(x_n) .
	\end{split}\]
	Note that the vector $\vec w = \paran{ w_1, \ldots, w_N}$ is independent of $x$ and is computed during the creation of the kernel matrix. While the computation is straightforward, a numerical issue arises when the test point $x$ is far from the data set $\SetDef{x_n}{n\in 1, \ldots, N}$. In that case, the vector $\paran{ k_{Gauss, \epsilon} (x,x_n)
	}_{n=1}^{N}$ have all entries below the machine precision. As a result the constant $\vec{\rho}_r(x)$ is also below machine precision. The occurrence of such near-zero numbers could lead to several computational anomalies. We use a nearest point approximation for evaluating the kernel sum in such cases. This also leads to false fixed points emerging at locations far from the dataset. We discuss this issue more in our concluding section.
	
	\paragraph{Summary of results} Figure \ref{fig:orbits} represents some typical SDE orbits that were provided as input to the algorithms. The results of the reconstruction have been displayed in Figures \ref{fig:Hopf}, \ref{fig:l63}, and \ref{fig:L96_N5} respectively. The expected trend that is displayed is that data created out of higher noise levels lead to larger errors in reconstruction. While the plots present the pointwise comparisons, the root mean squared errors are summarized in Table \ref{tab:results}. The trend can also be explained by the convergence analysis of the conditional expectation scheme.
	
	An important test for drift reconstruction algorithms is their efficacy in reconstructing trajectories.   To this end we have compared the trajectories generated by the true drift, with trajectories generated by the simulated drift. The trajectories show excellent match for the periodic Hopf oscillator (Figure \ref{fig:Hopf_recon}), and a close match for the Lorenz 63 model (Figure \ref{fig:L63_recon}). The Lorenz-63 model is a chaotic system, and the presence of positive Lyapunov exponents lead to an exponential divergence of trajectories. The gradual divergence of the trajectories in Figure \ref{fig:L63_recon} is indicative of a good reconstruction. It is also notable that the reconstructed orbit, even though divergent, remains topologically close to the true attractor $X$. The true Lorenz attractor is an isolated invariant set, with an attracting neighborhood. The topological proximity of the reconstruction indicates that the reconstructed drift preserves the topological confinement of a neighborhood of $X$.
	
	The 5-dimensional Lorenz-96 model presents a much steeper challenge, as evident from the trajectory comparison plots in Figures \ref{fig:L96_recon_a} and \ref{fig:L96_recon_b} respectively. As explained in Section \ref{sec:cnvrg} and Theorem \ref{thm:2}, the rate of convergence depends on the number of data samples $N$ as well as the standard deviation $\sigma$ of the noise. The higher the dimension $d$ of the phase-space, the higher $\sigma$ will be. This means a larger number of data samples are required to accurately sample the phase space. This is the well known and inevitable curse of dimensionality, given a fixed error bound $\delta$, the number of samples required to accurate sample the distribution scales as $\delta^{1/N}$. 
	We end this article by proposing a technique that can overcome this challenge on-parametrically, for certain drifts having redundancy in their structure.
	
	\section{Sparsity} \label{sec:sparse}
	
	The mechanism by which Algorithm \ref{algo:1} approximates the conditional expectations can be broadly described as computing averages with respect to various empirical measures. Empirical measure is the measure which is directly available to the algorithm. It is the discrete probability measure whose supports are the data points. The key principle behind the technique is that the empirical measure provides a weak approximation of the unknown measure $\mu$ referred to in Lemma \ref{lem:dj0k3}. In our next assumption, we make a precise statement on how each component of the drift is bound to a no more than $m$ other coordinates. This number $m$ is much smaller compared to the dimension $d$. Moreover, the functional dependence is also limited in type, as described  formally below :
	
	\begin{Assumption} \label{A:3}
		Given the SDE \eqref{eqn:sde}, there is 
		\begin{enumerate} [(i)]
			\item an integer $m<<d$ called the \emph{dependency dimension};
			\item a finite collection of functions $\tilde{V}^{l} : \real^m \to \real$, for $1\leq l \leq L$;
			\item and for each $i\in \braces{1, \ldots, d}$ a collection $\text{Input}(i)$ of $m$ coordinates $\paran{ k_{i,1} , \ldots, k_{i,m} } \subset \braces{1, \ldots, d}$ and an index $\ell(i) \in 1,\ldots L$;
		\end{enumerate}
		such that
		\begin{equation} \label{eqn:sparse_drift}
			V \paran{x}_i = \tilde{V}^{\ell(i)} \paran{ x_{k_{i,1}}, \ldots, x_{k_{i,m}} } , \quad \forall 1\leq i\leq d.
		\end{equation}
	\end{Assumption}
	
	Assumption \ref{A:3} says that effectively, each component of the drift $V$ is a repeated instance of a lower dimensional functional unit $\tilde{V}^{(l)}$. Such drifts occur in many natural complex phenomenon, such as natural complex dynamical systems \cite[e.g.]{chen2025critical}, networks of low dimensional dynamical systems \cite[e.g.]{yang2022controlling} and sub-grid modelling \cite[e.g.]{guan2023learning, Boutros2023onsager}. We now derive the precise nature of this expectation. For each $1\leq i\leq d$, define projection maps
	\[ \pi^i : \real^d \to \real^{m}, \quad x \mapsto \paran{ x_{k_{i,1}}, \ldots, x_{k_{i,m}} } .\]
	These projection maps create the following commutation with the coordinate projections $\SetDef{ \proj_i }{1\leq i \leq d}$, as shown below : 
	\begin{equation} \label{eqn:do9r}
		\begin{tikzcd} [column sep = large]
			\real^d \arrow[r, "\pi^i"] \arrow[d, "V"'] & \real^m \arrow[d, "\tilde{V}^{\ell(i)}"] \\
			\real^d \arrow[r, "\proj_i"'] & \real
		\end{tikzcd} , \quad \forall 1\leq i \leq d .
	\end{equation}
	Taking $\phi = \proj_i$ , we get the following equality with the aid of the limiting relation \eqref{eqn:deriv_exp:1} and commutation \eqref{eqn:do9r} :
	\[ \frac{d}{dt} \mathbb{E} \paran{ X_i( t, x_0) | x_0 } = \paran{ \proj_i V }(x_0) = \tilde{V}^{\ell(i)} \circ \pi^i (x_0) , \quad \forall 1\leq i \leq d .  \]
	
	Recall the function $f$ from \eqref{eqn:def:f}. Then similar to \eqref{eqn:deriv_exp:4} we can write 
	\begin{equation} \label{eqn:dhkmd0}
		V_i = \tilde{V}^{\ell(i)} = \mathbb{E}_{\mu} \paran{ \proj_i \circ f \,\rvert\, \pi_i } , \quad \forall 1\leq i \leq d.
	\end{equation}
	Equation \eqref{eqn:dhkmd0} says that each component of the drift $\tilde{V}$ can be recovered as a conditional expectation involving measurements of only a part of the coordinates of the full $d$-dimensional state space. It is also important to note that this equality holds for each $i$. This means that \eqref{eqn:dhkmd0} provides two advantages from a data-handling perspective. Firstly, even if $d$ is large, one does not need simultaneous measurements of each of the $d$ coordinates at every time step. Secondly, the number of coordinates or variables involved in the calculation is $m$, which could be much smaller than $d$. 
	
	To correctly utilize the strength of \eqref{eqn:dhkmd0}, we impose some restrictions on the data being used. A \emph{causally complete snapshot} corresponding to an initial condition $x_0$ and a coordinate $i\in 1,\ldots, d$ is an $(4+m)$-vector :
	\[
	\paran{ X^{(i)} \paran{ 3\Delta t, x_0 }, X^{(i)} \paran{ 2\Delta t, x_0 }, X^{(i)} \paran{ \Delta t, x_0 }, X^{(i)} \paran{ t, x_0 }, \pi^i \paran{ x_0 } } ,
	\]
	Each initial condition $x_0$ and choice of coordinate $i$ will lead to a different instantiation of this vector. The initial conditions need not be selected from a single path, they only need to be equidistributed w.r.t. $\mu_{stat}$. 
	
	The conditional expectation scheme \eqref{eqn:dhkmd0} can be numerically realized in a manner exactly analogous to Algorithm \ref{algo:2}. The convergence rates undergo an analysis similar to Section \ref{sec:cnvrg}. The scheme in \eqref{eqn:dhkmd0} allows the components of the drift to be mined from high dimensional data, in which trajectories may be short, and partially missing in data. A numerical implementation of \eqref{eqn:dhkmd0} to real-world high dimensional data is an important and promising project. These numerical investigations will be conducted in our forthcoming work.
	
	\section{Conclusions} \label{sec:conclus}
	
	Thus we have demonstrated using theory as well as numerical examples that the drift component of an SDE can be extracted from data, using principles from Probability theory and the theory of Kernel integral operators. The drift can be expressed as a conditional expectation involving random variables, on a certain probability space. This conditional expectation in turn can be estimated as a linear regression problem using techniques from kernel integral operators. The algorithmic procedure that is created is completely data-driven. It receives as inputs snapshots of trajectories of the unknown SDE, and the only requirement is that the initial condition in each snapshot is equidistributed with respect to some stationary measure of the SDE.  The numerical method worked well for the chaotic and quasiperiodic systems of dimensions less than 6 that we investigated. The performance deteriorates as the dimension increases.
	
	\paragraph{Choice of parameters} The three main tuning parameters in our algorithm are the Ridge regression $\theta$ and kernel bandwidth parameters $\epsilon_1, \epsilon_2, \epsilon_3$. Among these, the most important is the smoothing parameter $\epsilon_3$ which determines the effective radius of integration. The smaller $\epsilon_3$ is, the higher the number $N$ of total samples has to be. Thus a larger smoothing radius $\epsilon_3$ leads to a faster convergence with $N$. On the other hand, a larger $\epsilon_3$ and $N$ also increases the condition number of the Markov matrix $\Matrix{P}$, which could adversely affect the accuracy of the solution to the linear least squares problem \eqref{eqn:scheme:3}. There is no recipe for tuning these parameters that would work uniformly well in all applications. 
	
	\paragraph{Choice of kernel} Note that Algorithm~\ref{algo:1} does not specify the kernel, and any RKHS-like kernel such as the diffusion kernel would be sufficient. The convergence laid down in Lemma \ref{lem:dj0k3} also does not depend on the choice of kernel. However when finite data is used, the choice of kernel becomes important. The question of an optimal kernel is an important and unsolved question in machine learning \cite[e.g.]{narayan2021optimal, baraniuk1993signal, crammer2002kernel}. Our framework for estimating drift provides a more objective setting for this question.
	
	\paragraph{Out-of-sample extensions} The technique that we presented provides a robust estimation technique for the drift. The reason it is well suited to a general nonlinear and data-driven technique is the use of kernel methods. The downside of kernel methods is the difficulty in extending the function beyond the compactly supported dataset. Thus our technique provides a good phase-portrait, but it may not be reliable as a simulator or alternative model to the true system if the latter is chaotic and high dimensional. The issue of out of sample extensions \cite[e.g.]{VuralGuillemot2016out, PanEtAl2016out} and topological containment \cite[e.g.]{kaczynski2016towards, Mrozek2021creating} of simulation models is separate from the preliminary task of estimating the drift. This is yet another important area of development. 
	
	The present article reconstructs the drift as a sum of kernel sections, i.e., as a function. Kernel methods are just one of the many ways of representing a function. The next goal is to obtain a different representation that is more suitable to retaining the topology, but which still aims to achieve the identity \eqref{eqn:deriv_exp:4}. 
	A recent work by the author \cite{Das2024zero} presents a reconstruction strategy that preserves the topology of the targeted attractor, at the cost of replacing the deterministic dynamics law by a Markov process with arbitrarily low noise. The continuous time analog of this strategy is the idea of combinatorial drifts by Mrozek et. al. \cite{Mrozek1996multivalued, MischaikowMrozek1995isolating}. Another promising technique directed towards preserving the attractor topology is the matrix based approach of \cite{du2024computation}. A fourth option is the recent idea of representing vector fields as 1-forms \cite{DGGS2024sec}. A combination of any of these topologically motivated techniques, with the statistical methods presented in this article, offers an interesting avenue for further research.
	
	\paragraph{Acknowledgments} The author is grateful to the referees and editors for the detailed comments and reflections on the paper. They helped to improve the paper significantly.
	
	\bibliographystyle{unsrt}
	\bibliography{\Path References,ref}

\begin{thebibliography}{10}

\bibitem{Arnold1974stoch}
L.~Arnold.
\newblock Wiley Interscience, 1974.

\bibitem{philipp2023error}
F.~Philipp, M.~Schaller, K.~Worthmann, S.~Peitz, and F.~N{\"u}ske.
\newblock Error bounds for kernel-based approximations of the {K}oopman
  operator, 2023.

\bibitem{Doob1953book}
M.~Doob.
\newblock {\em Stochastic processes}.
\newblock Wiley, 1953.

\bibitem{HuangEtAl2015integral}
W.~Huang, M.~Ji, Z.~Liu, and Y.~Yi.
\newblock Integral identity and measure estimates for stationary
  {F}okker--{P}lanck equations.
\newblock {\em Ann. Probab.}, 43:1712--1730, 2015.

\bibitem{harris1956existence}
T.~Harris.
\newblock The existence of stationary measures for certain {M}arkov processes.
\newblock In {\em Proceedings of the Third Berkeley Symposium on Mathematical
  Statistics and Probability}, volume~2, pages 113--124, 1956.

\bibitem{HuangEtAl2018concn}
W.~Huang, M.~Ji, Z.~Liu, and Y.~Yi.
\newblock Concentration and limit behaviors of stationary measures.
\newblock {\em Physica D:}, 369:1--17, 2018.

\bibitem{ShaliziKontorovich2010}
C.~Shalizi and A.~Kontorovich.
\newblock {\em Almost none of the theory of stochastic processes}.
\newblock Lecture Notes, 2010.

\bibitem{Stanton1997nonpar}
R.~Stanton.
\newblock A nonparametric model of term structure dynamics and the market price
  of interest rate risk.
\newblock {\em J. Finance}, 52(5):1973--2002, 1997.

\bibitem{gouesbet1991reconstruction}
G~Gouesbet.
\newblock Reconstruction of the vector fields of continuous dynamical systems
  from numerical scalar time series.
\newblock {\em Physical Review A}, 43(10):5321, 1991.

\bibitem{tsutsumi2022constructing}
N.~Tsutsumi, K.~Nakai, and Y.~Saiki.
\newblock Constructing differential equations using only a scalar time-series
  about continuous time chaotic dynamics.
\newblock {\em Chaos}, 32(9), 2022.

\bibitem{gouesbet1994global}
G~Gouesbet and Christophe Letellier.
\newblock Global vector-field reconstruction by using a multivariate polynomial
  l 2 approximation on nets.
\newblock {\em Phys. Rev. E}, 49(6):4955, 1994.

\bibitem{cao2011robust}
J~Cao, L~Wang, and J~Xu.
\newblock Robust estimation for ordinary differential equation models.
\newblock {\em Biometrics}, 67(4):1305--1313, 2011.

\bibitem{ait2010operator}
Y.~A{\"\i}t-Sahalia, L.~Hansen, and J.~Scheinkman.
\newblock Operator methods for continuous-time {M}arkov processes.
\newblock In {\em Handbook of financial econometrics: tools and techniques},
  pages 1--66. Elsevier, 2010.

\bibitem{Garcia2017nonpar}
C.~Garcia et~al.
\newblock Nonparametric estimation of stochastic differential equations with
  sparse gaussian processes.
\newblock {\em Phys. Rev. E}, 96(2):022104, 2017.

\bibitem{Devlin2019opt}
J.~Devlin, D.~Husmeier, and J.~Mackenzie.
\newblock Optimal estimation of drift and diffusion coefficients in the
  presence of static localization error.
\newblock {\em Phys. Rev. E}, 100(2):022134, 2019.

\bibitem{darcy2023one}
M.~Darcy et~al.
\newblock One-shot learning of stochastic differential equations with data
  adapted kernels.
\newblock {\em Physica D}, 444:133583, 2023.

\bibitem{Batz2018approx}
P.~Batz, A.~Ruttor, and M.~Opper.
\newblock Approximate {B}ayes learning of stochastic differential equations.
\newblock {\em Phys. Rev. E}, 98(2):022109, 2018.

\bibitem{RuttorEtAl2013approximate}
A.~Ruttor, P.~Batz, and M.~Opper.
\newblock Approximate gaussian process inference for the drift function in
  stochastic differential equations.
\newblock {\em Adv. Neural Inf. Processing Sys.}, 26, 2013.

\bibitem{ella2024nonparametric}
E.~Ella-Mintsa.
\newblock Nonparametric estimation of the diffusion coefficient from iid sde
  paths.
\newblock {\em Statistical Inference for Stochastic Processes}, 27(3):585--640,
  2024.

\bibitem{davis2022est}
W.~Davis and B.~Buffett.
\newblock Estimation of drift and diffusion functions from unevenly sampled
  time-series data.
\newblock {\em Phys. Rev. E}, 106(1):014140, 2022.

\bibitem{ye2024nonlinear}
F.~Ye, S.~Yang, and M.~Maggioni.
\newblock Nonlinear model reduction for slow--fast stochastic systems near
  unknown invariant manifolds.
\newblock {\em J. Nonlinear Sci.}, 34(1):22, 2024.

\bibitem{Wang2020magnus}
Z.~Wang et~al.
\newblock The magnus expansion for stochastic differential equations.
\newblock {\em J. Nonlinear Sci.}, 30:419--447, 2020.

\bibitem{yin2022backstepping}
Xin X.~Yin and Qichun Zhang.
\newblock Backstepping-based state estimation for a class of stochastic
  nonlinear systems.
\newblock {\em Complex Engineering Systems}, 2(1):N--A, 2022.

\bibitem{Das2023conditional}
S.~Das.
\newblock Conditional expectation using compactification operators.
\newblock {\em Appl. Comput. Harmon. Anal.}, 71:101638, 2024.

\bibitem{DasGiannakis_delay_2019}
S.~Das and D.~Giannakis.
\newblock Delay-coordinate maps and the spectra of {K}oopman operators.
\newblock {\em J. Stat. Phys.}, 175:1107–1145, 2019.

\bibitem{DasGiannakis_RKHS_2018}
S.~Das and D.~Giannakis.
\newblock Koopman spectra in reproducing kernel {H}ilbert spaces.
\newblock {\em Appl. Comput. Harmon. Anal.}, 49:573--607, 2020.

\bibitem{DasMustAgar2023_qpd}
S.~Mustavee, S.~Das, and S.~Agarwal.
\newblock Data-driven discovery of quasiperiodically driven dynamics.
\newblock {\em Non. Dyn.}, Data-driven Nonlinear and Stochastic Dynamics with
  Control, 2024.

\bibitem{DasEtAl2023traffic}
S.~Das, S.~Mustavee, S.~Agarwal, and S.~Hassan.
\newblock Koopman-theoretic modeling of quasiperiodically driven systems:
  Example of signalized traffic corridor.
\newblock {\em IEE Trans. SMC Sys.}, 53:4466--4476, 2023.

\bibitem{DGJ_compactV_2018}
D.~Giannakis, S.~Das, and J.~Slawinska.
\newblock Reproducing kernel {H}ilbert space compactification of unitary
  evolution groups.
\newblock {\em Appl. Comput. Harmon. Anal.}, 54:75--136, 2021.

\bibitem{MarshallOlkin1960Cheby}
A.~Marshall and I.~Olkin.
\newblock Multivariate chebyshev inequalities.
\newblock {\em Annals Math. Stat.}, pages 1001--1014, 1960.

\bibitem{LorenzEmanuel1998opti}
E.~Lorenz and K.~Emanuel.
\newblock Optimal sites for supplementary weather observations: Simulation with
  a small model.
\newblock {\em J. Atmos. Sci.}, 55(3):399--414, 1998.

\bibitem{MarshallCoifman2019}
N.~Marshall and R.~Coifman.
\newblock Manifold learning with bi-stochastic kernels.
\newblock {\em IMA J. Appl. Math.}, 84(3):455--482, 2019.

\bibitem{WormellReich2021}
C.~Wormell and S.~Reich.
\newblock Spectral convergence of diffusion maps: Improved error bounds and an
  alternative normalization.
\newblock {\em SIAM J. Numer. Analy.}, 59(3):1687--1734, 2021.

\bibitem{CoifmanLafon2006}
R.~Coifman and S.~Lafon.
\newblock Diffusion maps.
\newblock {\em Appl. Comput. Harmon. Anal.}, 21:5–30, 2006.

\bibitem{HeinEtAl2005}
M.~Hein, JY. Audibert, and U.~Von Luxburg.
\newblock From graphs to manifolds--weak and strong pointwise consistency of
  graph {L}aplacians.
\newblock In {\em International Conference on Computational Learning Theory},
  pages 470--485. Springer, 2005.

\bibitem{VaughnBerryAntil2019}
R.~Vaughn, T.~Berry, and H.~Antil.
\newblock Diffusion maps for embedded manifolds with boundary with applications
  to pdes.
\newblock {\em Appl. Comput. Harmonic Anal.}, 68:101593, 2024.

\bibitem{BerryDas2023learning}
T.~Berry and S.~Das.
\newblock Learning theory for dynamical systems.
\newblock {\em SIAM J. Appl. Dyn.}, 22:2082 -- 2122, 2023.

\bibitem{BerryDas2024review}
T.~Berry and S.~Das.
\newblock Limits of learning dynamical systems.
\newblock {\em SIAM review}, 16, 2025.

\bibitem{chen2025critical}
D.~Chen et~al.
\newblock Critical nodes identification in complex networks: a survey.
\newblock {\em Complex Engineering Systems}, 5, 2025.

\bibitem{yang2022controlling}
CL. Yang and C.~Suh.
\newblock On controlling dynamic complex networks.
\newblock {\em Physica D}, 441:133499, 2022.

\bibitem{guan2023learning}
Y.~Guan et~al.
\newblock Learning physics-constrained subgrid-scale closures in the small-data
  regime for stable and accurate les.
\newblock {\em Physica D}, 443:133568, 2023.

\bibitem{Boutros2023onsager}
D.~Boutros and E.~Titi.
\newblock Onsager’s conjecture for subgrid scale $\alpha$-models of
  turbulence.
\newblock {\em Physica D}, 443:133553, 2023.

\bibitem{narayan2021optimal}
A.~Narayan, L.~Yan, and T.~Zhou.
\newblock Optimal design for kernel interpolation: Applications to uncertainty
  quantification.
\newblock {\em J. Comput. Phys.}, 430:110094, 2021.

\bibitem{baraniuk1993signal}
R.~Baraniuk and D.~Jones.
\newblock A signal-dependent time-frequency representation: optimal kernel
  design.
\newblock {\em IEEE Trans. Signal Process.}, 41(4):1589--1602, 1993.

\bibitem{crammer2002kernel}
K.~Crammer, J.~Keshet, and Y.~Singer.
\newblock Kernel design using boosting.
\newblock {\em Adv. Neural Inf Process. Sys.}, 15, 2002.

\bibitem{VuralGuillemot2016out}
E.~Vural and C.~Guillemot.
\newblock Out-of-sample generalizations for supervised manifold learning for
  classification.
\newblock {\em IEEE Trans. Image Processing}, 25(3):1410--1424, 2016.

\bibitem{PanEtAl2016out}
B.~Pan et~al.
\newblock Out-of-sample extensions for non-parametric kernel methods.
\newblock {\em IEEE Trans. neural netw. learning sys.}, 28(2):334--345, 2016.

\bibitem{kaczynski2016towards}
T.~Kaczynski, M.~Mrozek, and T.~Wanner.
\newblock Towards aformal tie between combinatorial and classical vector field
  dynamics.
\newblock {\em J. Computational Dynamics}, 3(1):17--50, 2016.

\bibitem{Mrozek2021creating}
M.~Mrozek and T.~Wanner.
\newblock Creating semiflows on simplicial complexes from combinatorial vector
  fields.
\newblock {\em J. Dif. Eq.}, 304:375--434, 2021.

\bibitem{Das2024zero}
S.~Das.
\newblock Reconstructing dynamical systems as zero-noise limits, 2024.

\bibitem{Mrozek1996multivalued}
M.~Mrozek.
\newblock Topological invariants, mulitvalued maps and computer assisted proofs
  in dynamics.
\newblock {\em Computers Math. Appl.}, 32(4):83--104, 1996.

\bibitem{MischaikowMrozek1995isolating}
K.~Mischaikow and M.~Mrozek.
\newblock Isolating neighborhoods and chaos.
\newblock {\em Japan J. Industrial Appl. Math.}, 12:205--236, 1995.

\bibitem{du2024computation}
Q.~Du and H.~Yang.
\newblock Computation of robust positively invariant set based on direct
  data-driven approach.
\newblock {\em Complex Engineering Systems}, 4(4):N--A, 2024.

\bibitem{DGGS2024sec}
S.~Das, D.~Giannakis, Y.~Gu, and J.~Slawinska.
\newblock Learning dynamical systems with the spectral exterior calculus, 2025.

\end{thebibliography}
\end{document}